\documentclass[a4paper]{amsart}
\usepackage{amsmath,amssymb,mathtools}
\usepackage{algorithm, algorithmic}

\newtheorem{theorem}{Theorem}[section]
\newtheorem{lemma}[theorem]{Lemma}
\newtheorem{corollary}[theorem]{Corollary}
\newtheorem{proposition}[theorem]{Proposition}

\theoremstyle{definition}
\newtheorem{definition}[theorem]{Definition}
\newtheorem{example}[theorem]{Example}
\newtheorem{remark}[theorem]{Remark}

\newtheorem{assumption}[theorem]{Assumption}

\DeclareMathOperator{\Ann}{Ann}

\DeclareMathOperator{\End}{End}
\DeclareMathOperator{\Sym}{Sym}
\DeclareMathOperator{\Cen}{Z}
\DeclareMathOperator{\Spec}{Spec}

\newcommand\LM{\mathop{\rm LM}\nolimits}
\newcommand\LT{\mathop{\rm LT}\nolimits}
\newcommand\LC{\mathop{\rm LC}\nolimits}
\newcommand\lcm{\mathop{\rm lcm}\nolimits}
\newcommand\Mat{\mathop{\rm Mat}\nolimits}
\newcommand\tsum{\textstyle\sum\limits}

\let\phi=\varphi
\let\rho=\varrho
\newcommand\ZZ{{\mathbb{Z}}}
\newcommand\QQ{{\mathbb{Q}}}
\newcommand\FF{{\mathbb{F}}}

{}

\title{Decomposing Finite $\mathbb{Z}$-Algebras}

\author{Martin Kreuzer}
\address[Martin Kreuzer]{Fakult\"{a}t f\"{u}r Informatik und Mathematik \\
Universit\"{a}t Passau, D-94030 Passau, Germany}
\email{martin.kreuzer@uni-passau.de}

\author{Alexei Miasnikov}
\address[Alexei Miasnikov]{Department of Mathematical Sciences,
Stevens Institute of Technology, 1 Castle Point Terrace,
Hoboken, NJ 07030, USA}
\email{amiasnikov@gmail.com}

\author{Florian Walsh}
\address[Florian Walsh]{Fakult\"{a}t f\"{u}r Informatik und Mathematik \\
Universit\"{a}t Passau, D-94030 Passau, Germany}
\email{florian.walsh@uni-passau.de}

\begin{document}

\begin{abstract}
For a finite $\ZZ$-algebra~$R$, i.e., for a ring which is not necessarily
associative or unitary, but whose additive group is finitely generated,
we construct a decomposition of $R/\Ann(R)$ into directly indecomposable factors under 
weak hypotheses. The method is based on constructing and decomposing a ring of scalars~$S$,
and then lifting the decomposition of~$S$ to the bilinear map given by the 
multiplication of~$R$, and finally to~$R/\Ann(R)$. All steps of the construction 
are given as explicit algorithms
and it is shown that the entire procedure has a probabilistic polynomial time 
complexity in the bit size of the input, except for the possible need to calculate
the prime factorization of one integer.
In particular, in the case when $\Ann(R) = 0$, these algorithms compute direct
decompositions of $R$ into directly indecomposable factors.
\end{abstract}



\maketitle


\section{Introduction}

In the early 1980s, A. Miasnikov and V. N. Remeslennikov studied
the elementary theory of finite-dimensional algebras (not necessarily associative,
not necessarily unitary) over regularly separable fields (see~\cite{MR1}, \cite{MR2}), 
making use of Mal'cev's correspondence (see~\cite{Mal}). Later, in the early 1990s,
A. Miasnikov extended and improved these results to finite-dimensional algebras
over arbitrary fields (cf.~\cite{Mya2}), and even to bilinear maps satisfying some
weak hypotheses (cf.~\cite{Mya1}). In the last decade this topic has received
renewed interest and many improvements and generalizations have been found 
(see for instance \cite{MS1} and \cite{GMO1}). 

For a finite $\ZZ$-algebra~$R$, i.e., a ring which is not necessarily associative 
or unitary, but whose additive group is finitely generated, 
a central step is the construction of a maximal ring of scalars 
which is a commutative, unitary ring. Then one uses the theory of commutative rings to
write the maximal ring of scalars as a direct product of indecomposable rings, and finally
this decomposition is lifted to a decomposition of~$R/\Ann(R)$.

The goal of this paper is to make all steps of this process explicit
and algorithmic, and to keep carefully track of the complexity of all operations
involved in the process. Already the case $R = \ZZ/n\ZZ$, where $n=pq$ is the product 
of two large prime numbers $p$ and~$q$, shows that the complexity of the entire 
construction cannot be purely polynomial time. However, as we shall see,
everything can be computed in probabilistic polynomial time plus one integer
factorization.

Let us discuss the results of the paper in more detail.
In Section~2 we present an algorithm to compute a presentation
via generators and relations of the maximal ring of scalars of a bilinear map.
Let $N_1, N_2, M$ be abelian groups and assume that $f:\; N_1 \times N_2 \longrightarrow M$ 
is a bilinear map.
The maximal ring of scalars $\mathfrak{S}(f)$ of~$f$ is a subring of 
the ring $\End(N_1) \times \End(N_2)$
and can be defined by equations over $N_1 \times N_2$. After fixing presentations 
of the groups $N_1$, $N_2$ and~$M$, its defining equations can be found by solving 
linear diophantine equations. This means that the maximal ring of scalars can be computed 
in polynomial time. After making all steps in the construction of $\mathfrak{S}(f)$,
we present the resulting algorithm in Proposition~\ref{mrs_endos}, and then we turn
its result into a $\ZZ$-algebra presentation of $\mathfrak{S}(f)$ in Proposition~\ref{mrs_pres}.

Next, in Section~3, we apply this technique to compute a presentation of the maximal
ring of scalars of a finite $\ZZ$-algebra~$R$. We can calculate the maximal ring of scalars
for the bilinear map
$$
f:\; R/\Ann_\lambda(R) \;\times\; R/\Ann_\rho(R) \;\longrightarrow\; R^2
$$ 
and derive the maximal ring of scalars $\mathfrak{S}(R)$ of~$R$ from~$\mathfrak{S}(f)$
(see Proposition~\ref{prop:linear_condition}).

The next step, namely the product decomposition of a finite unitary and commutative $\ZZ$-algebra,
is taken in Section~4. For this task, we need to solve the problem of computing the 
primitive idempotents of an explicitly presented commutative finite $\ZZ$-algebra
$S = \ZZ[x_1,\dots,x_n]/I$, where $I$ is an ideal in $P=\ZZ[x_1,\dots,x_n]$.
To this end, we use a simplified version of the algorithm given in~\cite{PSS}
for the computation of a primary decomposition of~$I$ (see Algorithm~\ref{decomp} 
and Proposition~\ref{prop:decomp}). Then it only remains to compute the connected
components of $\Spec(S)$ (see Algorithm~\ref{alg_connected}) and thus the primitive idempotents
of~$S$ (see Proposition~\ref{alg_idempotents}). For a presentation of~$S$ like the one we
have for a maximal ring of scalars, all algorithms in this section can be performed
in probabilistic polynomial time in the bit size of the input presentation, except for
one integer prime factorization (see~\cite{KW}).

In Section~5 we combine the results of Section~2 and Section~4. More precisely,
given a bilinear map $f:\; N_1\times N_2 \longrightarrow M$ as above, we compute
a presentation of its maximal ring of scalars $\mathfrak{S}(f)$ and then find
the primitive idempotents of~$\mathfrak{S}(f)$. Using these, we decompose~$f$
into a direct product of directly indecomposable bilinear maps 
(see Corollary~\ref{cor:bilinear_decomp}).

The topic of the last section is to lift this decomposition from the multiplication
map of a finite $\ZZ$-algebra to a decomposition of the algebra itself.
In this case we obtain a decomposition of the residue class ring $R/\Ann(R)$
into a direct product of $\ZZ$-algebras (see Theorem~\ref{decompose_R}). The algorithm
for computing this decomposition has again probabilistic polynomial time complexity
plus one integer prime factorization. Using suitably constructed finite dimensional Lie rings,
we see that some decompositions of~$R/\Ann(R)$ may not be constructable by lifting
a decomposition of~$\mathfrak{S}(R)$ (see Example~\ref{ex:not_a_lifting}), 
and that the decomposition of~$R/\Ann(R)$ cannot be lifted to a decomposition
of~$R$ without further hypotheses (see Examples~\ref{ex:not_liftable1} and~\ref{ex:not_liftable2}).
Moreover, Proposition~\ref{prop:indecomp_factors} provides a criterion under which the factors 
in the decomposition given by Theorem~\ref{decompose_R} are directly indecomposable.

One of the main applications of these results is that if $\Ann(R) = 0$ in a finite $\ZZ$-algebra $R$, then the
algorithms designed in this paper compute direct decompositions of $R$ into directly indecomposable factors and
these computations are quite efficient.

Unless explicitly noted otherwise, we use the definitions and notation introduced
in~\cite{KR1} and~\cite{KR3}. The algorithms in this paper were implemented
by the third author in the computer algebra system {\tt ApCoCoA} (see~\cite{ApCoCoA}). 
The source code is available from the authors upon request.

\bigbreak
%
%

\section{The Maximal Ring of Scalars of a Bilinear Map}%
\label{sec:max_ring}

In this section we study the maximal ring of scalars of a bilinear map
between abelian groups. The existence proof follows the construction in~\cite{Mya1}.
Our main goal is to find an algorithm for calculating the maximal ring of scalars.

\begin{definition}
Let $N_1$, $N_2$ and $M$ be abelian groups and $f \colon N_1 \times N_2 \longrightarrow M$ a map.
\begin{enumerate}
\item[(a)] The map $f$ is called {\bf bilinear} if, for every $a\in N_1$ and $b\in N_2$,
the maps $f(a, \cdot)$ and $f(\cdot, b)$ are group homomorphisms.

\item[(b)] The map $f$ is called {\bf non-degenerate} if $f(a,b)=0$ for all $b\in N_2$
implies $a=0$, and if $f(a,b)=0$ for all $a\in N_1$ implies $b=0$.

\item[(c)] If the additive group generated by $f(N_1,N_2)$ is~$M$ then the map~$f$
is called {\bf full}.
\end{enumerate}
\end{definition}

The central object of study in this section is defined as follows.

\begin{definition}\label{def_mrs}
Let $N_1, N_2,M$ be abelian groups and $f: N_1 \times N_2 \longrightarrow M$ a
bilinear map.
\begin{enumerate}
\item[(a)] A commutative ring $R$ is called a {\bf ring of scalars} of~$f$
if~$N_1$, $N_2$ and~$M$ are faithful $R$-modules, and if $f(ra, b) = f(a, rb) = r\, f(a,b)$
for all $r\in R$, $a\in N_1$ and $b \in N_2$.

\item[(b)] A ring of scalars~$R$ of~$f$ is called {\bf maximal} if every other
ring of scalars $R'$ of~$f$ can be embedded as a subring of $R$.
\end{enumerate}
\end{definition}

A ring of scalars can be embedded into a larger ring as follows.

\begin{remark}\label{rmk_embedding}
Let $N_1, N_2 ,M$ be abelian groups and $f:\; N_1\times N_2 \longrightarrow M$
a bilinear map.
\begin{enumerate}
\item[(a)] By viewing~$N_1$, $N_2$ and~$M$ as $\ZZ$-modules,
    we have the rings of endomorphisms $\End(N_1)$, $\End(N_2)$ and $\End(M)$.

\item[(b)] Given a ring of scalars~$R$ of~$f$, we associate
to every element $r\in R$ the {\bf multiplication endomorphisms}
$\phi_{i,r}:\; N_i \longrightarrow N_i$ given by $\phi_{i,r}(a) = ra$ for $a\in N_i$, and
$\psi_r:\; M \longrightarrow M$ given by $\psi_r(c)=rc$ for $c\in M$.

\item[(c)] Since $N_1$, $N_2$ and~$M$ are faithful $R$-modules,
the maps $\Phi_i:\; R \longrightarrow \End(N_i)$ given by
$\Phi_i(r)=\phi_{i,r}$ and $\Psi:\; R \longrightarrow \End(M)$
given by $\Psi(r)=\psi_r$ are injective ring homomorphisms.
\end{enumerate}
\end{remark}

In the following we let $N_1,N_2,M$ be abelian groups,
let $f: N_1\times N_2 \longrightarrow M$ be a bilinear map,
and let~$R$ be a ring of scalars of~$f$. Using the maps~$\Phi_i$
in part~(c) of the preceding remark, we identify~$R$ with its images
in~$\End(N_1)$ and $\End(N_2)$, and diagonally in $\End(N_1) \times \End(N_2)$.
Our next goal is to show that a maximal ring of scalars of~$f$ exists.

\begin{definition}
The elements of the sets
\begin{align*}
\Sym(f) &:= \{(\phi_1, \phi_2) \in \End(N_1) \times \End(N_2) \mid f(\phi_1(a), b) = f(a,\phi_2(b)) \\
        &\qquad\text{ for all } a \in N_1, b \in N_2 \} \; \text{ and } \\
\Cen(\Sym(f)) &:= \{\phi \in \Sym(f) \mid \phi \psi = \psi \phi \text{ for all } \psi \in \Sym(f)\}
\end{align*}
are called the \textbf{symmetric}, respectively \textbf{central symmetric}, endomorphisms of $N_1 \times N_2$
with respect to~$f$.
\end{definition}

The set of central symmetric endomorphisms of~$N_1 \times N_2$ has the following properties.

\begin{proposition}\label{Sym}
Let $N_1, N_2 ,M$ be abelian groups, and let $f: N_1\times N_2 \longrightarrow M$
be a bilinear map.
\begin{enumerate}
\item[(a)] The set $\Cen(\Sym(f))$ is a commutative unitary subring of $\End(N_1) \times \End(N_2)$.

\item[(b)] If $f$ is non-degenerate then every ring of scalars $R$ of $f$
can be embedded as subring of $\Cen(\Sym(f))$.
\end{enumerate}
\end{proposition}

\begin{proof}
To prove~(a), we need to show that $\Cen(\Sym(f))$ is closed under addition and
multiplication. For all $\phi = (\phi_1, \phi_2),\psi = (\psi_1, \psi_2) \in \Cen(\Sym(f))$
and all $a\in N_1$ and $b \in N_2$, we have
\begin{align*}
    f((\phi_1+\psi_1)(a),b) = f(\phi_1(a),b)+f(\psi_1(a),b) &= f(a,\phi_2(b)) + f(a,\psi_2(b)) \\
                                                            &= f(a,(\phi_2+\psi_2)(b))
\end{align*}
and $f(\phi_1\psi_1(a), b) = f(\psi_1(a),\phi_2(b)) = f(a,\psi_2\phi_2(b)) = f(a,\phi_2 \psi_2(b))$.
Consequently, we get $\phi+\psi \in \Cen(\Sym(f))$ and $\phi\psi \in \Cen(\Sym(f))$, as claimed.

To show~(b), we let~$R$ be a ring of scalars of~$f$. Notice that, by definition,
we have $R \subseteq \Sym(f)$. It remains to show that $R \subseteq \Cen(\Sym(f))$.
For $\phi=(\phi_1, \phi_2)\in R$ and $\psi = (\psi_1, \psi_2) \in \Sym(f)$ and all
$a\in N_1$ and $b \in N_2$, we have
$$
f(\phi_1\psi_1(a),b) = \phi\cdot f(\psi_1(a),b) = \phi\cdot f(a,\psi_2(b)) =
f(\phi_1(a),\psi_2(b)) = f(\psi_1\phi_1(a),b)
$$
Thus we obtain $f((\phi_1\psi_1 - \psi_1\phi_1)(a),b)=0$ for all $a\in N_1$ and $b \in N_2$.
Since $f$ is non-degenerate, it follows that $(\phi_1\psi_1 - \psi_1\phi_1)(a)=0$ for all $a\in N_1$.
Analogously we obtain $(\phi_2\psi_2 - \psi_2\phi_2)(b)=0$ for all $b\in N_2$, and hence
$\phi\psi = \psi\phi$. This proves $R\subseteq \Cen(\Sym(f))$.
\end{proof}

Now we are ready to prove the existence of a maximal ring of scalars of~$f$.
The following theorem was first shown in~\cite{Mya1}, Thm.~1, using a
different setting and notation.

\begin{theorem}\label{max_ring_existence}
Let $N_1, N_2 ,M$ be abelian groups, and let $f: N_1 \times N_2 \longrightarrow M$
be a non-degenerate and full bilinear map.
Then there exists a unique maximal ring of scalars of~$f$.
It is denoted by~$\mathfrak{S}(f)$ and satisfies
\begin{align*}
\mathfrak{S}(f) \;=\; \{ &(\phi_1, \phi_2) \in \Cen(\Sym(f)) \mid {\tsum}_i f(\phi_1(a_i),b_i) =
{\tsum}_j f(\phi_1(a'_j),b'_j)\\
&\text{for }a_i,a'_j \in N_1 \text{ and } b_i,b'_j\in N_2
\text{ with }{\tsum}_i f(a_i,b_i) = {\tsum}_j f(a'_j,b'_j) \in M \}.
\end{align*}
\end{theorem}

\begin{proof}
The ring $\Cen(\Sym(f))$ is commutative and acts faithfully on~$N_1$ and $N_2$.
By Proposition~\ref{Sym}, we know that every ring of scalars is a subring of~$\Cen(\Sym(f))$.
Since the map~$f$ is full, every element $c\in M$ can be represented
in the form $c = \sum_i f(a_i,b_i)$ with $a_i \in N_1$ and $b_i \in N_2$.

Let $R$ be the set of all maps $(\phi_1, \phi_2) \in \Cen(\Sym(f))$ such that
$\sum_i f(a_i,b_i) = \sum_j f(a'_j,b'_j)$ implies
$\sum_i f(\phi_1(a_i),b_i) = \sum_j f(\phi_1(a'_j),b'_j)$
for all $a'_j \in N_1$ and $b'_j\in N_2$ with $c=\sum_j f(a'_j,b'_j)$.
Using the bilinearity of $f$, we see that for $\phi = (\phi_1, \phi_2),\psi = (\psi_1, \psi_2) \in R$
and $a_i,b_i,a'_j,b'_j,c$ as above, we have
\begin{align*}
{\tsum}_i f((\phi_1+\psi_1)(a_i),b_i) &= {\tsum}_i f(\phi_1(a_i),b_i) + {\tsum}_i f(\psi_1(a_i),b_i) =\\
&= {\tsum}_j f(\phi_1(a'_j),b'_j) + {\tsum}_j f(\psi_1(a'_j),b'_j) \\
&= {\tsum}_j f((\phi_1+\psi_1)(a'_j),b'_j)
\end{align*}
and therefore $\phi+\psi\in R$. Furthermore, ${\tsum}_i f(\psi_1(a_i),b_i) =
{\tsum}_j f(\psi_1(a'_j),b'_j)$ implies ${\tsum}_i f(\phi_1\psi_1(a_i),b_i) =
{\tsum}_j f(\phi_1\psi_1(a'_j),b'_j)$
Hence the set $R$ is closed under addition and multiplication,
and therefore a subring of~$\Cen(\Sym(f))$.

By definition, the action $\phi(c) := \sum_i f(\phi(a_i),b_i)$ on
$c=\sum_i f(a_i,b_i) \in M$ is well-defined for all $\phi \in R$.
Since the action of~$R$ on~$N_1$ and $N_2$ is faithful
and~$f$ is non-degenerate, the action of~$R$ on~$M$ is faithful as well.
Thus~$R$ is a ring of scalars of~$f$.

It remains to show that~$R$ is a maximal ring of scalars and that it is unique.
Let $\tilde{R}$ be a ring of scalars of~$f$. Using Remark~\ref{rmk_embedding},
we identify~$\tilde{R}$ with its image in~$\End(N_1) \times \End(N_2)$. Then,
by Proposition~\ref{Sym}, we have $\tilde{R} \subseteq \Cen(\Sym(f))$. Finally,
since~$M$ is a faithful $\tilde{R}$-module, the ring~$\tilde{R}$ has to satisfy
the defining condition of~$R$, and thus is a subring of~$R$. This shows that~$R$
is a maximal ring of scalars of~$f$. As every other maximal ring of scalars
embeds into~$R$, it is the unique ring with this property.
\end{proof}

In the following we denote the additive group of a ring $R$ by~$R^+$.
In the next step we describe the maximal ring of scalars $\mathfrak{S}(f)$ of~$f$
as a set of solutions of a system of homogeneous linear equations over~$\ZZ$.
Since it is a subring of the ring $\End(N_1) \times \End(N_2)$, the
first goal is to describe the additive groups of $\End(N_1)$ and $\End(N_2)$ as
the set of solutions of such a linear system.

In the following we assume that $N$ is a finitely generated abelian group and that
a presentation of~$N$ is given by
$$
N \;=\; \ZZ e_1 \oplus\cdots\oplus \ZZ e_n \;/\;
\big\langle\, {\tsum_{i=1}^n} r_{i1}\, e_i,\; \dots,\; {\tsum_{i=1}^n} r_{i\nu}\,e_i \big\rangle
$$
with $r_{ij} \in \ZZ$ for $i=1,\dots,n$ and $j=1,\dots,\nu$.
For $i=1,\dots,n$, we denote the residue class
of~$e_i$ in~$N$ by~$a_i$. Then $\{a_1,\dots,a_n\}$ is a system of
generators of~$N$. We shall describe $\End(N)$ by solving a
homogeneous linear system over~$N$ in the following sense.

\begin{remark}{\bf (Solving Systems of Linear Equations over~$N$)}%
\label{solve_over_N}\\
Let $N$ be a finitely generated abelian group with a presentation
as above, let $v_1,\dots,v_p\in N$, and let
\begin{equation*}
\tag{i}\label{eqn_over_N} x_1 v_1 + \cdots + x_p v_p \;=\; 0
\end{equation*}
be a homogeneous linear equation over~$N$ in the indeterminates
$x_1,\dots,x_p$. For $i=1,\dots,p$, we write
$v_i = c_{i1}a_1 + \cdots + c_{in} a_n$ with $c_{ij}\in\ZZ$.

Then a tuple $(s_1,\dots,s_p) \in \ZZ^p$ is a solution of~(\ref{eqn_over_N})
if and only if there exist $s_{p+1},\dots,s_{p+\nu}\in\ZZ$ such that
$(s_1,\dots,s_{p+\nu})$ is a solution of the system of~$n$ homogeneous
linear equations
\begin{equation*}
\tag{ii}\label{eqn_over_Z}  \tsum_{i=1}^p x_i (c_{i1}e_1 + \cdots + c_{in} e_n)
+ \tsum_{j=1}^\nu x_{p+j} (r_{1j}e_1 + \cdots + r_{nj} e_n) \;=\; 0
\end{equation*}
where $x_{p+1},\dots,x_{p+\nu}$ are further indeterminates.

Hence we may find a $\ZZ$-basis of the set of solutions of equation~(\ref{eqn_over_N})
by solving the system~(\ref{eqn_over_Z}), projecting the generators of the
solution space to their first~$p$ coordinates, and interreducing the result.
\end{remark}

Now we are ready to describe $\End(N)$ as follows.

\begin{proposition}\label{end_gen}
Let $N$ be a finitely generated abelian group with a presentation
as above. Consider the following system of homogeneous linear equations over~$N$
in the indeterminates $x_{ij}$, where $i,j \in \{1,\dots,n\}$:
\begin{equation*}
\tag{iii}\label{eqEnd}
r_{1i}\; {\textstyle\sum\limits_{j=1}^n} x_{1j} a_j \;+\; \cdots
\;+\; r_{ni}\; {\textstyle\sum\limits_{j=1}^n} x_{nj} a_j \;=\; 0
\end{equation*}
for $i=1,\dots,\nu$.
Then the elements of~\/$\End(N)$ correspond to the solutions of {\rm (\ref{eqEnd})}.
More precisely, a solution $(c_{ij}) \in (\ZZ^n)^n$ of {\rm (\ref{eqEnd})}
corresponds to the endomorphism $\phi \in \End(N)$ which satisfies
$\phi(a_i) = \sum_{j=1}^n c_{ij} a_j$ for $i=1,\dots,n$.
\end{proposition}

\begin{proof}
To every $\phi \in \End(N)$, we associate the $n$-tuple
$(b_1, \dots, b_n)\in N^n$ with $b_i=\phi(a_i)$
for $i=1,\dots,n$. An arbitrary tuple $(b_1, \dots, b_n) \in N^n$
is associated to an endomorphism of~$N$
if and only if it satisfies the relations of~$N$, i.e., iff
$\sum_{i=1}^n r_{ij} b_i = 0$ for $j=1,\dots,\nu$.
We write each component $b_i$ as a linear combination of the
generators $\{a_1,\dots,a_n\}$ of~$N$
and get $b_i = \sum_{j=1}^n c_{ij} a_j$ with coefficients $c_{ij}\in \ZZ$.
Thus the tuple $(b_1,\dots,b_n)$ is associated to an element
$\phi \in \End(N)$ if and only if the tuple of coefficients $(c_{ij})$
is a solution of~(\ref{eqEnd}).
\end{proof}

Having described the elements of $\End(N_1) \times \End(N_2)$ as solutions of a system of homogeneous
linear equations, we now consider $\Sym(f)$, $\Cen(\Sym(f))$, and $\mathfrak{S}(f)$.
Let $N_1$ be generated by $\{a_1, \dots, a_n\}$ and let $N_2$ be generated by $\{a_{n+1}, \dots, a_{n'}\}$.
Since these sets depend also on the abelian group~$M$, we need to fix
a presentation
$$
M \;=\; \ZZ e_1 \oplus \cdots\oplus \ZZ e_m \;/\;
\big\langle\, {\tsum_{i=1}^m} r'_{i1}\, e_i,\; \dots,\; {\tsum_{i=1}^m} r'_{i\mu}\,e_i \big\rangle
$$
with $r'_{ij}\in\ZZ$ for $i=1,\dots,m$ and $j=1,\dots,\mu$.
The residue class of $e_j$ in~$M$ is denoted by~$b_j$ for $j=1,\dots,m$.
Then the set $\{b_1,\dots,b_m\}$ is a system of generators of~$M$.

Furthermore, we may describe~$f$ explicitly as follows.
For $i=1,\dots,n$ and $j=n+1, \dots, n'$, we write $f(a_i,a_j)=\gamma_{ij1}b_1 +\cdots +
\gamma_{ijm}b_m$ with $\gamma_{ij1},\dots,\gamma_{ijm}\in \ZZ$.
Notice that we can now solve the systems of homogeneous linear equations
over~$M$ appearing in the following proposition in analogy to
the method explained in Remark~\ref{solve_over_N}. In the following we assume that elements
$\phi = (\phi_1, \phi_2) \in\End(N_1) \times \End(N_2)$ are defined under the bijection
given in Proposition~\ref{end_gen} by tuples $(d_{ij})$ in $(\ZZ^{n'})^{n'}$ with
$\phi_1(a_i) = \sum_{j=1}^n d_{ij} a_j$ for $i=1, \dots, n$ and
$\phi_2(a_i) = \sum_{j=n+1}^{n'} d_{ij} a_j$ for $i=n+1, \dots, n'$.

\begin{proposition}\label{rings_eq}
Let $(d_{ij})$ be a tuple in~$(\ZZ^{n'})^{n'}$ which is associated to an element
$\phi\in\End(N_1) \times \End(N_2)$ as above. Then the elements of $\Sym(f)$, $Z(\Sym(f))$,
and $\mathfrak{S}(f)$ correspond to the solutions of the following systems of
homogeneous linear equations.
\begin{enumerate}
\item[(a)] The element~$\phi$ is contained in $\Sym(f)$ if and only if $(d_{ij}) \in (\ZZ^{n'})^{n'}$
    is a solution of the following system of homogeneous linear equations over~$M$:
$$
{\tsum_{i=1}^n} \; x_{k i} f(a_i, a_\ell) \;-\;
{\tsum_{j=n+1}^{n'}} \; x_{\ell j} f(a_k, a_j) \;=\; 0
$$
for $k\in \{1,\dots,n\}$ and $\ell\in \{n+1,\dots,n'\}$

\item[(b)] Let $(\delta_{ij}^{\,(1)}),\dots, (\delta_{ij}^{\,(p)})$
be tuples in $(\ZZ^{n'})^{n'}$ such that the corresponding endomorphisms
in~$\End(N_1) \times \End(N_2)$ generate the additive group $\Sym(f)$.
Then an element $\phi$ of $\Sym(f)$ is contained in $Z(\Sym(f))$ if and only if
the tuple $(d_{ij})$ is a solution of the following system of homogeneous linear equations over~$M$:
$$
{\tsum_{i=1}^n}{\tsum_{j=n+1}^{n'}} \; x_{k i}\, \delta_{\ell j}^{\,(\pi)} f(a_i,a_j)
\;-\; {\tsum_{i=1}^n}{\tsum_{j=n+1}^{n'}} \; x_{\ell j}\, \delta_{k i}^{\,(\pi)} f(a_i,a_j) \;=\; 0
$$
for $k\in \{1,\dots,n\}$, $\ell \in \{n+1, \dots, n'\}$ and $\pi\in \{1,\dots,p\}$.

\item[(c)] Let $(z_{ij}^{\,(1)}), \dots, (z_{ij}^{\,(q)})$ be tuples in
$(\ZZ^n)^{n'-n}$ which form a $\ZZ$-basis of the solution space of the
following homogeneous linear equation over~$M$:
$$
{\tsum_{i=1}^n} {\tsum_{j=1}^{n'-n}} \; x_{ij}\, f(a_i,a_{j+n}) \;=\; 0
$$
Then an element $\phi$ of $Z(\Sym(f))$ given by a tuple~$(d_{ij}) \in (\ZZ^{n'})^{n'}$ is contained
in~$\mathfrak{S}(f)$ if and only if $(d_{ij})_{i=1, \dots,n}^{j=1,\dots,n}$ in $(\mathbb{Z}^n)^n$
is a solution of the following system of homogeneous linear equations over~$M$:
$$
{\tsum_{i,\ell=1}^n} {\tsum_{j=1}^{n'-n}}\; x_{i\ell}\; z_{ij}^{\,(k)}\;
f(a_\ell, a_{j+n}) \;=\; 0
$$
for $k\in\{1,\dots,q\}$.
\end{enumerate}
\end{proposition}

\begin{proof}
A version of these claims was shown in~\cite{Mya1}. For the convenience of the reader,
and since it is the basis for any implementation,
we provide a proof using our notation.

To prove~(a) we note that, for a tuple $(d_{ij})\in (\ZZ^{n'})^{n'}$,
we have $\sum_{i=1}^n d_{k i} f(a_i, a_\ell) =
f(\phi_1(a_k),a_\ell)$ and $\sum_{j=n+1}^{n'}
d_{\ell j} f(a_k, a_j) = f(a_k, \phi_2(a_\ell))$.
Hence the tuple $(d_{ij})$ solves the system if and only if~$\phi$
is symmetric.

Next we show~(b). For $\pi=1,\dots,p$, let $\psi^\pi = (\psi_1^\pi, \psi_2^\pi)$
be the tuple of endomorphisms corresponding to $(\delta_{ij}^{\,(\pi)})$
By assumption, $\{ \psi^1,\dots,\psi^p\}$ generates $\Sym(f)$.
For every tuple $(d_{ij}) \in (\ZZ^{n'})^{n'}$, we have
\begin{align*}
    \tsum_{i=1}^n \tsum_{j=n+1}^{n'} d_{k i} \delta_{\ell j}^{\,(\pi)} f(a_i,a_j) &=
    f(\phi_1(a_k),\psi_2^\pi(b_\ell)) \text{ and}\\
    \tsum_{i=1}^n \tsum_{j=n+1}^{n'} d_{\ell j} \delta_{k i}^{\,(\pi)} f(a_i,a_j) &=
    f(\psi_1^\pi(a_k),\phi_2(a_\ell)).
\end{align*}
Therefore this tuple
solves the system if and only if we have
$f(\phi_1(a_k),\rho_2(b_\ell)) = f(\rho_1(a_k),\phi_2(b_\ell))$
for every $\rho = (\rho_1, \rho_2) \in\Sym(f)$. Since~$\rho$ is symmetric, this is
equivalent to
$$
f(\rho_1\phi_1(a_k),b_\ell) = f(\phi_1(a_k),\rho_2(b_\ell)) = f(\rho_1(a_k),
\phi_2(b_\ell)) = f(\phi_1\rho_1(a_k),b_\ell).
$$
Now we use the facts that~$f$ is
non-degenerate and that $\{a_{n+1},\dots,a_{n'}\}$ generates~$N_2$ to conclude
that the latter condition is equivalent to $\rho_1\phi_1=\phi_1\rho_1$.
Since $\phi$ is symmetric, we analogously obtain $\rho_2\phi_2=\phi_2\rho_2$,
which proves $\phi\in Z(\Sym(f))$.

Finally, we prove~(c). By writing the elements $a_i,b_i,a'_j,b'_j$
in Theorem~\ref{max_ring_existence} in terms of the generators of~$N_1$ and $N_2$,
it is easy to check that~$\phi$ is contained in~$\mathfrak{S}(f)$ if and only if
$$\tsum_{i=1}^n \tsum_{j=1}^{n'-n} \gamma_{ij} f(a_i,a_{j+n})=0 \quad \text{implies} \quad
\tsum_{i=1}^n \tsum_{j=1}^{n'-n} \gamma_{ij} f(\phi_1(a_i),a_{j+n})=0$$
for all $(\gamma_{ij}) \in (\ZZ^n)^{n'-n}$. By construction, the tuples $(z_{ij}^{\,(k)})$
with $k\in \{1,\dots,q\}$ form a $\ZZ$-basis of the solution space
of the first equation. For every $(d_{ij})\in (\ZZ^{n'})^{n'}$, the tuple
$(d_{ij})_{i=1, \dots, n}^{j=1, \dots, n}$ satisfies
$$
\tsum_{i,\ell=1}^n \tsum_{j=n+1}^{n'} \; d_{i\ell}\; z_{ij}^{\,(k)}\;
f(a_\ell, a_j) = \tsum_{i=1}^n \tsum_{j=n+1}^{n'} z_{ij}^{\,(k)} \, f(\phi_1(a_i),b_j).
$$
Hence we have
$\sum_{i=1}^n \sum_{j=n+1}^{n'} z_{ij}^{\,(k)} \, f(\phi_1(a_i),a_j)=0$ for $k=1,\dots,q$,
if and only if this tuple solves the second system in~(c)
and this is equivalent to $\phi\in \mathfrak{S}(f)$.
\end{proof}

Part~(c) of the preceding proposition allows us to compute
tuples representing a system of generators of~$\mathfrak{S}(f)$.
It remains to determine the relations between these generators.
This can be done as follows.

\begin{proposition}\label{relations}
Let $(d_{ij}^{\,(1)}),\dots, (d_{ij}^{\,(\ell)})$ be tuples in $(\ZZ^{n'})^{n'}$ which
represent elements $\phi_1,\dots,\phi_\ell \in \mathfrak{S}(f)$ with $\phi_i = (\phi_{i1}, \phi_{i2})$,
and assume that $\{\phi_1,\dots,\phi_\ell\}$ generates the $\ZZ$-module $\mathfrak{S}(f)^+$.
Then a tuple of integers $(\gamma_1,\dots,\gamma_\ell)\in \ZZ^\ell$
represents a relation $\sum_{i=1}^\ell \gamma_i \phi_i =0$
if and only if it is a solution of the system of homogeneous linear equations over $N_1 \times N_2$ given by
$$
{\tsum_{i=1}^\ell}{\tsum_{k=1}^n} \; x_i\, d_{jk}^{\,(i)}\, a_k \;=\; 0
$$
for $j=1,\dots,n$ and
$$
{\tsum_{i=1}^\ell}{\tsum_{k=n+1}^{n'}} \; x_i\, d_{jk}^{\,(i)}\, a_k \;=\; 0
$$
for $j=n+1,\dots,n'$.
\end{proposition}

\begin{proof}
Clearly, we have $\sum_{i=1}^\ell \gamma_i \phi_i=0$ if and only if
we have $\sum_{i=1}^\ell \gamma_i \phi_{i1}(a_j)=0$ for $j=1,\dots,n$
and $\sum_{i=1}^\ell \gamma_i \phi_{i2}(a_j)=0$ for $j=n+1,\dots,n'$.
Now it suffices to use the equalities $\phi_{i1}(a_j) = \sum_{k=1}^n
d_{j k}^{\,(i)}\, a_k$ and $\phi_{i2}(a_j) = \sum_{k=n+1}^{n'} d_{j k}^{\,(i)}\, a_k$
\end{proof}

In summary, to compute a presentation of the maximal ring of scalars,
we need to solve the systems of linear equations over~$N_1$ and $N_2$ given in
Propositions~\ref{end_gen} and~\ref{relations}, and the systems over $M$ given in \ref{rings_eq}.
As we have noted in Remark~\ref{solve_over_N}, this amounts
to solving various systems of homogeneous linear equations over~$\ZZ$.
Let us mention one way to perform this well-known task.

\begin{remark}\label{linalgZcomplexity}
Let $A \in \Mat_{m,n}(\ZZ)$ be an integer matrix which defines a system of
homogeneous linear equations. Since the set of solutions of this system forms a
$\ZZ$-submodule of~$\ZZ^n$, it is free.
As for instance described in~\cite{La}, we can compute
a $\ZZ$-basis of the set of solutions of the system as follows.

For the matrix~$A$, there exist unimodular matrices $L \in \Mat_m(\ZZ)$
and $R \in \Mat_n(\ZZ)$ such that $LAR = D$ where $D$ is a diagonal
matrix called the {\it Smith normal form} of~$A$. Using these matrices, we
can describe all solutions of the linear system. Namely, let~$s$ be the rank of~$D$.
Then a $\ZZ$-basis of the set of solutions of the system is given
by $\{Re_{s+1},\dots,Re_n\}$. Hence it suffices to compute the Smith normal form
of~$A$ and the corresponding unimodular matrix~$R$.

In \cite{KaBa}, Kannan and Bachem give an algorithm for computing the Smith normal form
and the corresponding unimodular matrices in polynomial time and such that the number of digits
of all matrix entries are bounded by a polynomial. More detailed complexity analyses of
solving systems of linear equations over~$\ZZ$ are given by Storjohann in~\cite{St} 
and Wan in~\cite{Wan}.
\end{remark}

At this point we are ready to combine the results of this section
to give an algorithm which computes a system of generators
of the maximal ring of scalars of a bilinear map.

\begin{remark}\label{remark:alg_input}
As an input for all further calculations, we assume that we are given the following information.
\begin{enumerate}
\item[(a)] Presentations of~$N_1$ and $N_2$ of the form
\begin{align*}
    N_1 &= \ZZ e_1 \oplus \cdots \oplus \ZZ e_n  \mathbin{/}
    \langle \tsum_{i=1}^n r_{i1} e_i, \dots, \tsum_{i=1}^n r_{i\nu} e_i \rangle \quad \text{and} \\
    N_2 &= \ZZ e_1 \oplus \cdots \oplus \ZZ e_{n'-n}  \mathbin{/}
    \langle \tsum_{i=n+1}^{n'} r_{i1} e_{i-n}, \dots, \tsum_{i=n+1}^{n'} r_{i\nu'} e_{i-n} \rangle
\end{align*}
with $r_{ij}\in\ZZ$.
Denoting the residue class of~$e_i$ in~$N_1$ by~$a_i$ for $i=1,\dots,n$,
we get a system of generators $\{a_1,\dots,a_n\}$ for the $\ZZ$-module~$N_1$.
Analogously, denoting the residue class of~$e_i$ in~$N_2$ by~$a_{n+i}$ for $i=1,\dots,n'-n$,
we get a system of generators $\{a_{n+1},\dots,a_{n'}\}$ for the $\ZZ$-module~$N_2$.

\item[(b)] A presentation of~$M$ of the form
$$
M = \ZZ e_1 \oplus \cdots \oplus \ZZ e_m \mathbin{/}
\langle \tsum_{i=1}^m \rho_{i1} e_i, \dots, \tsum_{i=1}^m \rho_{i\mu} e_i \rangle
$$
with $\rho_{ij}\in\ZZ$. Denoting the residue class of~$e_i$ in~$M$ by~$c_i$ for $i=1,\dots,m$,
we get a system of generators $\{b_1,\dots,b_m\}$ for the $\ZZ$-module $M$.

\item[(c)] In order to encode the bilinear map~$f$, we write the image $f(a_i, a_j)$ as a
$\ZZ$-linear combination of the generators of~$M$ and obtain
$$
f(a_i, a_j) = \tsum_{k=1}^m s_{ijk} b_k
$$
with structure constants $s_{ijk}\in \ZZ$ for $i=1,\dots,n$, $j=n+1, \dots, n'$ and $k = 1, \dots, m$.
\end{enumerate}
\end{remark}

\begin{proposition}\label{mrs_endos}
Let $f \colon N_1 \times N_2 \longrightarrow M$ be a full and non-degenerate bilinear map.
Assume that we are given presentations as in Remark~\ref{remark:alg_input}.
Then the following steps define an algorithm which computes
elements $\phi_1,\dots,\phi_r$ in~$\End(N_1) \times \End(N_2)$ corresponding to a system of
$\ZZ$-module generators of~$\mathfrak{S}(f)$.
\begin{enumerate}
\item[(1)] Solve the system of homogeneous linear equations over~$N_1 \times N_2$ given by
\begin{equation*}
\tsum_{i,j=1}^n x_{ij} r_{ik} a_j   \;=\; 0
\end{equation*}
for $k\in\{1,\dots,\nu\}$ and
\begin{equation*}
    \tsum_{i,j=n+1}^{n'} x_{ij} r_{ik} a_j   \;=\; 0
\end{equation*}
for $k\in\{1,\dots,\nu'\}$, and let $L$ be its solution space.

\item[(2)] Solve the system of homogeneous linear equations over~$M$ given by
\begin{equation*}
\tsum_{i=1}^n x_{ki} \tsum_{\kappa=1}^m s_{i\ell\kappa} b_\kappa
- \tsum_{j=n+1}^{n'} x_{\ell j} \tsum_{\kappa=1}^m s_{kj\kappa} b_\kappa \;=\; 0
\eqno{(i)}
\end{equation*}
for $k\in \{1,\dots,n\}$ and $\ell \in \{n+1, \dots,n'\}$.

\item[(3)] Compute the intersection of $L$ and the solution space of the linear system~(i).
Let $((\delta_{ij}^{\,(1)}), \dots, (\delta_{ij}^{\,(p)})$ be tuples in $(\ZZ^{n'})^{n'}$
which form a $\ZZ$-basis of the intersection.

\item[(4)] Solve the system of homogeneous linear equations over~$M$ given by
\begin{equation*}
    \tsum_{i=1}^n \tsum_{j=n+1}^{n'} x_{ki} \delta_{\ell j}^{\, (\pi)} \tsum_{\kappa=1}^m s_{ij\kappa} b_\kappa
    \;-\; \tsum_{i=1}^n \tsum_{j=n+1}^{n'} x_{\ell j} \delta_{ki}^{\,(\pi)} \tsum_{\kappa=1}^m s_{ij\kappa} b_\kappa \;=\; 0
\eqno{(ii)}
\end{equation*}
where $k\in \{1,\dots,n\}$, $\ell \in \{n+1, \dots, n'\}$ and $\pi\in\{1,\dots,p\}$

\item[(5)] Solve the following homogeneous linear equation over~$M$:
\begin{equation*}
\tsum_{i=1}^n \tsum_{j=1}^{n'-n} x_{ij} \tsum_{\kappa=1}^m s_{ij\kappa} b_\kappa = 0
\eqno{(iii)}
\end{equation*}
Let $(z_{ij}^{\,(1)}),\dots, (z_{ij}^{\,(q)})$ be tuples in $(\ZZ^n)^{n'-n}$
which form a $\ZZ$-basis of the solution space.

\item[(6)] Solve the system of homogeneous linear equations over~$M$ given by
\begin{equation*}
    \tsum_{i,\ell=1}^n \tsum_{j=1}^{n'-n}  x_{i\ell} z_{ij}^{\,(k)} \tsum_{\kappa=1}^m s_{\ell j\kappa} b_\kappa
\;=\; 0
\eqno{(iv)}
\end{equation*}
for $k=1,\dots,q$ and let $\mathcal{L} \subseteq (\ZZ^n)^n$ be its solution space.

\item[(7)] Compute the intersection of $\mathcal{L} \oplus (\mathbb{Z}^{n'-n})^{n'-n}$ and the subspaces of
$(\ZZ^{n'})^{n'}$ obtained in Step~(3) and Step~(5). Let the tuples $(d_{ij}^{\,(1)}),\dots, (d_{ij}^{\,(r)})$
in $(\ZZ^{n'})^{n'}$ be a $\ZZ$-basis of the intersection.

\item[(8)] Return the elements $\phi_1,\dots,\phi_r$ in~$\End(N_1) \times \End(N_2)$
given by $\phi_k = (\phi_{k1}, \phi_{k2})$ with $\phi_{k1}(a_i)=\sum_{j=1}^n d_{ij}^{\,(k)} a_j$ and
$\phi_{k2}(a_i)=\sum_{j=n+1}^{n'} d_{ij}^{\,(k)} a_j$.
\end{enumerate}
\end{proposition}

\begin{proof}
By Proposition~\ref{end_gen}, the solutions of the system in Step~(1) correspond
to the tuples of endomorphisms in~$\End(N_1) \times \End(N_2)$. When we replace in~(i) the
sum $\sum_{\kappa=1}^m s_{ij\kappa} c_\kappa$ by~$f(a_i, a_j)$, we see that
it is exactly the system in Proposition~\ref{rings_eq}.a.
Therefore the intersection of $L$ and the solutions of~(i)
in Step~(3) corresponds to the elements of~$\Sym(f)$.
Next we note that system~(ii) in conjunction with the intersection in Step~(3) corresponds
by Proposition~\ref{rings_eq}.b to the system defining
the elements of $Z(\Sym(f))$.
Furthermore system~(iii) is equivalent to the hypothesis
in Proposition~\ref{rings_eq}.c and system~(iv) is the
conclusion of this statement. Together with Proposition~\ref{rings_eq}.c
we obtain that the intersection computed in Step~(7) represents
exactly the elements of~$\mathfrak{S}(f)$. Hence Step~(8) returns
the correct endomorphisms.
\end{proof}

The final step is to take the endomorphisms found in this
proposition and to compute a $\ZZ$-algebra presentation
of~$\mathfrak{S}(f)$ from them. This is achieved as follows.

\begin{proposition}\label{mrs_pres}
Let $f: N_1 \times N_2 \longrightarrow M$ be a full and non-degenerate
bilinear map and assume that we are given presentations as in Remark~\ref{remark:alg_input}.
Moreover, suppose we are given tuples $(d_{ij}^{\,(1)}),\dots, (d_{ij}^{\,(r)})$ in $(\ZZ^{n'})^{n'}$
which represent $\phi_1, \dots, \phi_r$ in~$\End(N_1) \times \End(N_2)$ corresponding to a system of
$\ZZ$-module generators of~$\mathfrak{S}(f)$.

Then the following steps define an algorithm which returns
tuples $(v_k)\in \ZZ^r$ and $(u_{\ell k}) \in \ZZ^{r}$
for $\ell=1,\dots,\rho$ and $(t_{ijk})\in \ZZ^r$ for $i,j=1,\dots,r$
such that we have a $\ZZ$-algebra presentation
$$
\mathfrak{S}(f) \,=\, \ZZ[y_1,\dots,y_r] \,/\,
\langle v_1y_1 + \dots + v_r y_r -1, {\tsum_{k=1}^r} u_{\ell k} y_k,\;
y_iy_j - {\tsum_{k=1}^r} t_{ijk} y_k \rangle
$$
for $i,j=1,\dots,r$ and $\ell=1,\dots,\rho$.
\begin{enumerate}
\item[(1)] For $i=1,\dots,r$, form the elements
$\delta_{ik} = \sum_{j=1}^n d_{kj}^{(i)} a_j \in N_1$ where $k \in \{1,\dots,n\}$ and
the elements $\delta_{ik} = \sum_{j=n+1}^{n'} d_{kj}^{(i)} b_j \in N_2$ where $k \in \{n+1,\dots,n'\}$.

\item[(2)] Solve the system of homogeneous linear equations over~$N_1 \times N_2$ given by
$$
\tsum_{i=1}^r x_i \delta_{ij} \;=\; 0
$$
for $j=1,\dots,n'$. Let $u_1, \dots, u_{\rho}$ be tuples in $\ZZ^{r}$ which
form a $\ZZ$-basis of the solution space.

\item[(3)] Solve the system of linear equations over~$N_1 \times N_2$ given by
$$
\tsum_{i=1}^r x_i \delta_{ik} - a_k \;=\; 0
$$
for $k=1, \dots, n'$. Let $(v_1, \dots, v_r)\in \ZZ^r$ be an
element of the solution space.

\item[(4)] For $i,j=1,\dots,r$, solve the system of homogeneous linear equations over $N_1 \times N_2$ given by
$$
\tsum_{\ell=1}^{n'} d_{k\ell}^{(j)} \delta_{i \ell} \;-\;
\tsum_{\kappa=1}^r x_{\kappa} \delta_{\kappa k} \;=\; 0
$$
for $k=1,\dots,n'$. For $i,j=1,\dots,r$, let $(t_{ij1}, \dots, t_{ijr})\in \ZZ^{r}$
be an element of the solution space.

\item[(5)] Return the tuples $(v_k)\in \ZZ^r$ and $(u_{\ell k}) \in \ZZ^r$
for $\ell\in\{1,\dots,\rho\}$ as well as $(t_{ijk}) \in \ZZ^r$ for
$i,j\in\{1,\dots,r\}$ and stop.

\end{enumerate}
\end{proposition}

\begin{proof}
In Step~(1) we form for $i=1, \dots, r$ the elements $\delta_{ik}$ in $N_1$ and $N_2$
such that we have $\phi_{i1}(a_k) = \delta_{ik}$ for $k=1,\dots,n$ and
$\phi_{i2}(a_k) = \delta_{ik}$ for $k=n+1,\dots,n'$.
Then the system of linear equations in Step~(2) corresponds to the system in
Proposition~\ref{relations}. Therefore Step~(2) computes the linear relations
of $\phi_1, \dots, \phi_r$.

Step~(3) computes a tuple $(v_1,\dots,v_r)\in \ZZ^r$ such that
$(v_1\phi_{11} + \cdots +v_r \phi_{r1})(a_k) = a_k$ for $k=1, \dots, n$ and
$(v_1\phi_{12} + \cdots +v_r \phi_{r2})(a_k) = a_k$ for $k=n+1, \dots, n'$.
Thus the linear combination $v_1\phi_1 + \cdots +v_r \phi_r$ is the identity map.
In order to finally get a $\ZZ$-algebra presentation of~$\mathfrak{S}(R)$,
it remains to express the products $\phi_i \phi_j$ using the system of generators
$\{\phi_1,\dots,\phi_r\}$. This is done in Step~(4), where we use
$$
\phi_{i1} \phi_{j1}(a_k) = \phi_{i1}(\tsum_{\ell=1}^n d_{k\ell}^{(j)} a_{\ell}) =
\tsum_{\ell=1}^n d_{k\ell}^{(j)} \phi_{i1}(a_{\ell}) = \tsum_{\ell=1}^n d_{k\ell}^{(j)} \delta_{i\ell}.
$$
Then for $i,j=1, \dots, r$ a solution $(t_{ij1}, \dots, t_{ijr})\in \ZZ^r$
of the system in Step~(4) satisfies
$\phi_{i1} \phi_{j1} (a_k) = (t_{ij1} \phi_{11} + \dots +  t_{ijr} \phi_{ri})(a_k)$ for $k=1, \dots, n$.
Analogously, a solution $(t_{ij1}, \dots, t_{ijr})\in \ZZ^r$ satisfies
$\phi_{i2} \phi_{j2} (a_k) = (t_{ij1} \phi_{12} + \dots +  t_{ijr} \phi_{r2})(a_k)$ for $k=n+1, \dots, n'$.
\end{proof}

Notice that all computations performed by the algorithms of
Proposition~\ref{mrs_endos} and Proposition~\ref{mrs_pres}
can be done in polynomial time with respect to the size
of the input and require only coefficients whose size is bounded
by a polynomial (see Remark~\ref{linalgZcomplexity}).

\bigbreak
%
%

\section{The Maximal Ring of Scalars of a Finite $\ZZ$-Algebra}%
\label{sec:max_ring_Zalg}

In this section we focus on our central topic, namely the
study of finite $\ZZ$-algebras.

\begin{definition}
Let~$(R,+)$ be an abelian group, i.e., a $\ZZ$-module.
\begin{enumerate}
\item[(a)] The group~$R$ is called a \textbf{$\ZZ$-algebra} if it is equipped with a
$\ZZ$-bilinear map $R \times R \longrightarrow R$. This bilinear map is called
the \textbf{multiplication} of~$R$ and is denoted by $(a,b) \mapsto ab$.

\item[(b)] A $\ZZ$-algebra is called \textbf{finite} if it is finitely generated 
as a $\ZZ$-module. Sometimes such rings are also called {\bf finite dimensional $\ZZ$-algebras},
or briefly {\bf FDZ-algebras}.

\end{enumerate}
\end{definition}

In other words, a $\ZZ$-algebra is a ring for which 
we do not require the existence of an identity element 
or associativity or commutativity of the multiplication. 
We denote the additive group of~$R$ by~$R^+$. 
It is a finitely generated abelian
group, and the multiplication in~$R$ gives rise to a bilinear map
$R^+ \times R^+ \longrightarrow R^2$, where~$R^2$ is the
subgroup of~$R^+$ generated by all products $a\cdot b$ with $a,b\in R$.
This bilinear map is clearly full, but it can be degenerate, and
therefore the results of Section~2 cannot be applied directly. In order to
construct a non-degenerate bilinear map, we proceed as follows.

\begin{definition}
Let~$R$ be a $\ZZ$-algebra. Then the sets
\begin{align*}
\Ann_\lambda(R) \;&=\; \{a\in R \mid a\, b = 0 \text{ for all } b \in R\},\\
\Ann_\rho(R) \;&=\; \{a\in R \mid b\, a = 0 \text{ for all } b \in R\},\\
\text{and}\qquad \Ann(R) \;&=\; \Ann_\lambda(R) \cap \Ann_\rho(R)
\end{align*}
are called the {\bf left}, {\bf right} and {\bf two-sided annihilator} of~$R$, respectively.
\end{definition}

Notice that $\Ann_\lambda(R)$ and $\Ann_\rho(R)$ are subgroups of~$R^+$ 
and that $\Ann(R)$ is a two-sided ideal in~$R$.

\begin{proposition}\label{nondegenerate_g}
Let $R$ be a $\ZZ$-algebra and let~$R^2$ be the subgroup of~$R^+$ generated
by the products $a\cdot b$, where $a,b\in R$. Then the map
$$
f_R:\; R^+/\Ann_\lambda(R) \times  R^+/\Ann_\rho(R) \;\longrightarrow\; R^2
$$
given by $f_R(a +\Ann_\lambda(R),b+\Ann_\rho(R)) = a\cdot b$ is non-degenerate and full.
\end{proposition}

\begin{proof}
Let us indicate the proof of the non-degeneracy of~$f_R$. Suppose that $a\in R$ is such that
its residue class~$\bar{a} = a+ \Ann_\lambda(R)$ satisfies $f_R(\bar{a},\bar{b}) = ab = 0$ for all
$\bar{b} \in R/\Ann_\rho(R)$. Then we have $ab=0$ for all $b\in R$, and thus $a\in \Ann_\lambda(R)$
implies $\bar{a}=0$. The non-degeneracy in the other argument follows in the same way.
\end{proof}

Now Theorem~\ref{max_ring_existence} shows that~$f_R$ has a
unique maximal ring of scalars $\mathfrak{S}(f_R)$. By definition of $\mathfrak{S}(f_R)$,
the $\ZZ$-modules $R^+/\Ann_\lambda(R)$, $R^+/\Ann_\rho(R)$,  and $R^2$ are faithful $\mathfrak{S}(f_R)$-modules.
But, in general, the action of $\mathfrak{S}(f_R)$ is not defined on all of~$R$. For our goal of computing a direct
decomposition of $R$ we need some additional properties of the ring $\mathfrak{S}(f_R)$.

\begin{definition}
Let $R$ be a $\ZZ$-algebra.
\begin{enumerate}
\item[(a)] A commutative ring~$S$ is called a \textbf{ring of scalars} of~$R$
if the following two conditions are satisfied.
\begin{enumerate}
\item[(i)] The ring~$S$ is a ring of scalars of the bilinear map
$$
f^*_R : R^+/\Ann(R) \times  R^+/\Ann(R) \;\longrightarrow\; R^2.
$$

\item[(ii)] The canonical $\ZZ$-linear map
$R^2 \longrightarrow R^+/\Ann(R)$ is $S$-linear.
\end{enumerate}

\item[(b)] A ring of scalars~$S$ of~$R$ is called \textbf{maximal}
if every other ring of scalars $S'$ of~$R$ can be embedded as a subring of~$S$.

\end{enumerate}
\end{definition}

To show the existence of the maximal ring of scalars of a $\ZZ$-algebra $R$ we can not
immediately apply Theorem~\ref{max_ring_existence}, since the bilinear map $f^*_R$ might be degenerate.
Instead we show that the maximal ring of scalars of $R$ is a subring of $\mathfrak{S}(f_R)$.

\begin{lemma}\label{ring_of_scalars_f_R}
    A ring of scalars of a $\ZZ$-algebra $R$ is a ring of scalars of the bilinear map $f_R$.
\end{lemma}
\begin{proof}
If $S$ is a ring of scalars of $R$ then $R/\Ann(R)$ and $R^2$ are faithful $S$-modules, the map
$f_R^*$ is $S$-bilinear and the map $R^2 \to R/\Ann(R)$ is $S$-linear. Observe that in this case the subgroups
$\Ann_\lambda(R)/\Ann(R)$ and $\Ann_\rho(R)/\Ann(R)$ are, in fact, $S$-submodules of $R/\Ann(R)$. To see this, let
$x \in \Ann_\lambda(R)$, let $\alpha \in S$, and suppose that $\alpha(x+\Ann(R)) = x' + \Ann(R)$ in $R/\Ann(R)$
for some $x' \in R$. For every $y \in R$ we get
$$
f_R^*(\alpha(x+\Ann(R)),y +\Ann(R)) = f_R^*(x'+\Ann(R),y+\Ann(R)) = x'y.
$$
On the other hand, we have
$$
f_R^*(\alpha(x+\Ann(R)),y +\Ann(R)) = \alpha(xy) = \alpha \cdot 0 = 0.
$$
Thus we obtain $x'y = 0$ for every $y \in R$, and hence $x' \in \Ann_\lambda (R)$. Therefore $\Ann_\lambda(R)/\Ann(R)$
and $\Ann_\rho(R)/\Ann(R)$ are $S$-submodules of $R/\Ann(R)$. This implies that the residue class modules
$$
(R/\Ann(R))/(\Ann_\lambda(R)/\Ann(R)) \quad \text{and} \quad (R/\Ann(R))/(\Ann_\rho(R)/\Ann(R))
$$
induce $S$-module structures on
$$
R/\Ann_\lambda(R) \cong (R/\Ann(R))/(\Ann_\lambda(R)/\Ann(R))
$$
and
$$
R/\Ann_\rho(R) \cong (R/\Ann(R))/(\Ann_\rho(R)/\Ann(R))
$$
in such a way that, in the notation above, if $\alpha(x +\Ann(R)) = x' + \Ann(R)$ in $R/\Ann(R)$ then
$\alpha(x +\Ann_\lambda(R)) = x' + \Ann_\lambda (R)$ in $R/\Ann_\lambda(R)$ and
$\alpha(x +\Ann_\rho(R)) = x' + \Ann_\rho (R)$ in $R/\Ann_\rho(R)$.
Note that $R^2$ is a faithful $S$-module, as mentioned above.

Now we are ready to show that the map $f_R$ is $S$-bilinear. For $y \in R$, we have
$$
f_R(\alpha(x +\Ann_\lambda(R)), y+\Ann_\rho(R)) = f_R(x' +\Ann_\lambda(R), y+\Ann_\rho(R)) = x'y.
$$
On the other hand,
$$
f_R^*(\alpha(x +\Ann(R)),y+\Ann(R)) = f_R^*(x' +\Ann(R), y+\Ann(R)) = x'y,
$$
so
$$
f_R(\alpha(x +\Ann_\lambda(R)), y+\Ann_\rho(R)) = f_R^*(\alpha(x +\Ann(R)),y+\Ann(R))
$$
for all $x, y \in R$. Hence, we deduce
\begin{align*}
    f_R(\alpha(x +\Ann_\lambda(R)), y+\Ann_\rho(R)) &=  f_R^*(\alpha(x +\Ann(R)),y+\Ann(R)) \\
                                                    &= \alpha f_R^*(x +\Ann(R),y+\Ann(R)) \\
                                                    &= \alpha (xy) = \alpha f_R(x+\Ann_\lambda(R),y+\Ann_\rho(R)).
\end{align*}

This shows that $f_R$ is $S$-linear in the first variable. A similar argument shows that $f_R$
is $S$-linear in the second variable.

Finally, it remains to show that $R/\Ann_\lambda(R)$ and $R/\Ann_\rho (R)$ are faithful $S$-modules.
Suppose that $\alpha \in S$ is such that for every $x \in R$ there exists an element $x' \in \Ann_\lambda(R)$
such that $\alpha(x+\Ann_\lambda(R)) = x'+ \Ann_\lambda(R)$. For every $y \in R$ we obtain
$$
f_R(\alpha(x+\Ann_\lambda(R)), y+\Ann_\lambda(R)) = \alpha(xy)
$$
and also
$$
f_R(\alpha(x+\Ann_\lambda(R)), y+\Ann_\lambda(R)) = f_R(x'+\Ann_\lambda(R), y+\Ann_\lambda(R)) = x'y = 0.
$$
Therefore we see that $\alpha(xy) = 0$ for all $x,y \in R$, so $\alpha R^2 = 0$. Since $R^2$ is a faithful $S$-module,
we conclude that $\alpha = 0$. This shows that $R/\Ann_\lambda(R)$ is a faithful $S$-module. A similar argument
proves that $R/\Ann_\rho(R)$ is a faithful $S$-module.
\end{proof}

The following theorem is a generalization of Proposition~8 in~\cite{Mya2}.

\begin{theorem}\label{max_ring_algebra_existence}
Let $R$ be a $\ZZ$-algebra, and let
$$\eta : R \rightarrow R^+/\Ann_\lambda(R) \times R^+/\Ann_\rho(R)$$ be the canonical diagonal group homomorphism.
Then there exists a unique maximal ring of scalars of $R$. It is the largest subring $S$
of $\mathfrak{S}(f_R)$ for which the following conditions hold.
\begin{itemize}
\item [(1)] The canonical group homomorphisms
$$
R^2 \rightarrow R^+/\Ann_\lambda(R) \quad \text{and} \quad R^2 \rightarrow R^+/\Ann_\rho(R)
$$
are $S$-linear, i.e.\ the restriction of $\eta$ to $R^2$ is $S$-linear.

\item [(2)] The image $\eta(R) \cong R^+/\Ann(R)$ is an $S$-submodule of the
$S$-module $$R^+/\Ann_\lambda(R) \times R^+/\Ann_\rho(R).$$
\end{itemize}
In the following we denote this ring~$S$ by $\mathfrak{S}(R)$.
\end{theorem}

\begin{proof}
Let $S$ be a ring of scalars of $R$. Then $S$ is a ring of scalars of $f_R$ by Lemma~\ref{ring_of_scalars_f_R}.
We start by proving that $S$ satisfies Conditions (1) and (2). Indeed, by the construction in
Lemma~\ref{ring_of_scalars_f_R}, the structure of an $S$-module on $R/\Ann_\lambda(R)$ and on
$R/\Ann_\rho(R)$ is induced by the structure of the $S$-modules
$$
(R/\Ann(R))/(\Ann_\lambda(R)/\Ann(R)) \quad \text{and} \quad (R/\Ann(R))/(\Ann_\rho(R)/\Ann(R)).
$$
Hence the group homomorphism
$$
R/\Ann(R) \to (R/\!\Ann(R))/(\Ann_\lambda(R)/\!\Ann(R)) \times (R/\!\Ann(R))/(\Ann_\rho(R)/\!\Ann(R))
$$
is $S$-linear. This proves that $S$ satisfies Condition~(2).

To show (1) observe that $R^2 \to R/\Ann(R)$ is $S$-linear, since $S$ is a ring of scalars of $R$.
Now the canonical group homomorphism
$R^2 \rightarrow R^+/\Ann_\lambda(R) \times  R^+/\Ann_\rho(R)$ is a composition of two $S$-linear homomorphisms
$$
R^2 \rightarrow R/\Ann(R)  \rightarrow R^+/\Ann_\lambda(R) \times  R^+/\Ann_\rho(R),
$$
hence it is $S$-linear, so (1) holds.

Conversely, let $S$ be a ring of scalars of $f_R$ which satisfies Conditions~(1) and~(2).
We now prove that $S$ is a ring of scalars of $f^*_R$. Since $S$ is a ring of scalars of $f_R$, the groups
$R/\Ann_\lambda(R)$, $R/\Ann_\rho(R)$, and $R^2$ are faithful $S$-modules. Condition~(2) ensures that
the embedding
$$
R/\Ann(R) \to R/\Ann_\lambda(R) \times R/\Ann_\rho(R)
$$
induces a faithful $S$-module structure on $R/\Ann(R)$ in such a way that, for $\alpha \in S$ and $x, x' \in R$
we have $\alpha(x+\Ann(R)) = x'+\Ann(R)$ if and only if $\alpha(x+\Ann_\lambda(R)) = x'+\Ann_\lambda(R)$ and
$\alpha(x+\Ann_\rho(R)) = x'+\Ann_\rho(R)$. Now we calculate
\begin{align*}
&f_R^*(\alpha(x+\Ann(R)), y+\Ann(R)) = f_R^*(x'+\Ann(R)), y+\Ann(R)) = x'y \\
&= f_R(x'+\Ann_\lambda(R)), y+\Ann_\rho(R)) = f_R(\alpha(x+\Ann_\lambda(R)), y+\Ann_\rho(R)) \\
&= \alpha f_R(x+\Ann_\lambda(R), y+\Ann_\rho(R)) = \alpha(xy) = \alpha f_R^*(x+\Ann(R), y+\Ann(R)).
\end{align*}
Hence $f_R^*$ is $S$-linear in the first variable. A similar argument shows that $f_R^*$ is $S$-linear in the second
variable. We proved that $S$ is a ring of scalars of $f_R^*$. It is left to show that the group homomorphism
$R^2 \to R/\Ann(R)$ is $S$-linear. But this follows from Condition~(1).

Finally, we need to show that there is a maximal ring $S$ of scalars of $f_R$ which satisfies Conditions~(1) and~(2).
Denote by $S$ the set of all elements $\alpha \in \mathfrak{S}(f_R)$
such that the maps $R^2 \rightarrow R^+/\Ann_\lambda(R)$ and $R^2 \rightarrow R^+/\Ann_\rho(R)$ are $\alpha$-linear and
the subset $\eta(R)$ of the $\mathfrak{S}(f_R)$-module $R^+/\Ann_\lambda(R) \times R^+/\Ann_\rho(R)$ is invariant under
the action of $\alpha$. It is easy to check that $S$ is a subring of $\mathfrak{S}(f_R)$. Hence it is the largest
subring of $\mathfrak{S}(f_R)$ satisfying Conditions~(1) and~(2). This proves the theorem.
\end{proof}

Let us now examine how to compute the maximal ring of scalars $\mathfrak{S}(R)$
of a finite $\ZZ$-algebra~$R$.
The first step is to apply the results of Section~2 and compute the maximal ring
of scalars of the bilinear map~$f_R$.

\begin{assumption}\label{remark:algebra_input}
For the remainder of this section, we assume that we are given the following information.
\begin{enumerate}
\item[(a)] A presentation of $R^+$ of the form
$$
R^+ \;=\; \ZZ e_1 \oplus\cdots\oplus \ZZ e_n \;/\;
\big\langle\, {\tsum_{i=1}^n} r_{i1}\, e_i,\; \dots,\;
{\tsum_{i=1}^n} r_{i\nu}\,e_i \big\rangle
$$
with $r_{ij}\in\ZZ$. Denoting the residue class of~$e_i$ in~$R^+$ by~$a_i$
for $i=1,\dots,n$, we get a system of generators $\{a_1,\dots,a_n\}$ for the
$\ZZ$-module $R^+$.
We denote the canonical image of $a_i$ in $R^+/\Ann_\lambda(R)$ by $b_i$ and the
canonical image of $a_i$ in $R^+/\Ann_\rho(R)$ by $c_i$.

\item[(b)] Furthermore, we assume that we are given equalities
$$
a_i a_j =  {\textstyle\sum\limits_{k=1}^n} \; s_{ijk} \, a_k
$$
with $s_{ijk}\in \ZZ$ for $i,j,k=1,\dots,n$.

\end{enumerate}
\end{assumption}

In this setting, we can calculate $\Ann_\lambda(R)$ and $\Ann_\rho(R)$ as follows.

\begin{remark}\label{present_Ann}
Let $R$ be a ring which is a finite $\ZZ$-module, and assume that~$R$
is presented as above. In order to find systems of $\ZZ$-module generators
for $\Ann_\lambda(R)$ and $\Ann_\rho(R)$, it suffices to solve
suitable systems of linear equations over~$R^+$ via Remark~\ref{solve_over_N}.

To compute a system of generators of~$\Ann_\lambda(R)$,
we have to find all $b\in R$ such that $b a_i=0$ for $i=1,\dots,n$.
Let us write $b=c_1a_1 +\cdots+ c_n a_n$ with $c_i\in \ZZ$.
Using the above equalities to replace $a_ia_j$, we see that
we have to solve the system of homogeneous linear equations
$$
x_1 {\tsum_{j=1}^n} s_{1ij} a_j \;+\; \cdots \;+\;
x_n {\tsum_{j=1}^n} s_{nij} a_j \;=\; 0
$$
where $i\in\{1,\dots,n\}$.

We can determine a presentation of $R^+/\Ann_\lambda(R)$. Let $\{b_1,\dots,b_\mu\}$ be a system
of generators of $\Ann_\lambda(R)$ and write $b_i=\sum_{j=1}^n d_{ij}a_j$ with $d_{ij}\in\ZZ$
for $i,j=1,\dots,n$. Then we get a presentation
$$
R^+/\Ann_\lambda(R) \cong \ZZ e_1 \oplus \cdots \oplus \ZZ e_n /
\big\langle {\tsum_{i=1}^n} r_{i1} e_i, \dots,
{\tsum_{i=1}^n} r_{i\nu}e_i, {\tsum_{j=1}^n} d_{1j}e_j, \dots,
{\tsum_{j=1}^n} d_{\mu j} e_j \big\rangle.
$$
A presentation of $R^+/\Ann_\rho(R)$ can be computed analogously.

\end{remark}

Having determined generators of the maximal ring of scalars of~$f_R$,
we can compute the subring corresponding to $\mathfrak{S}(R)$ as follows.

\begin{proposition}\label{prop:linear_condition}
Let $f_R$ be the bilinear map defined as above, and let $\phi = (\phi_1, \phi_2)$ be an element of its
maximal ring of scalars~$\mathfrak{S}(f_R)$. Assume that $\phi$ is given by a tuple
$(d_{ij}) \in (\ZZ^{2n})^{2n}$ defining endomorphisms $\phi_1$ of $R^+/\Ann_\lambda(R)$ by
$\phi_1(b_i) = \sum_{j=1}^n d_{ij} b_j$ and $\phi_2$ of $R^+/\Ann_\rho(R)$ by
$\phi_2(c_i) = \sum_{j=n+1}^{2n} d_{ij} c_{j-n}$
Then $\phi$ is an element of $\mathfrak{S}(R)$ if and only if the following conditions are satisfied.
\begin{enumerate}
    \item[(1)] The tuple $(d_{ij})$ is a solution of the system of homogeneous equations
        over $R^+/\Ann_\lambda(R) \times R^+/\Ann_\rho(R)$ given by
        \begin{align*}
            \tag{i} &\tsum_{\ell=1}^n \tsum_{k=1}^n s_{\ell j k} x_{i\ell}b_k
            - \tsum_{k=1}^n\tsum_{\ell=1}^n s_{ijk} x_{k\ell} b_{\ell} = 0 \quad \text{and} \\
            \tag{ii} &\tsum_{\ell=1}^n \tsum_{k=1}^n s_{\ell j k} x_{i+n, \ell+n}c_k
            - \tsum_{k=1}^n\tsum_{\ell=1}^n s_{ijk} x_{k+n, \ell+n} c_{\ell} = 0
        \end{align*}
        for $i,j \in \{1, \dots, n\}$.
    \item[(2)] There exist $\gamma_1, \dots, \gamma_n \in \mathbb{Z}$ such that the tuple $(d_{ij})$ is
        a solution of the system of equations over $R^+/\Ann_\lambda(R) \times R^+/\Ann_\rho(R)$ given by
        \begin{align*}
            \tsum_{j=1}^n x_{ij} b_j = \tsum_{j=1}^n \gamma_j b_j \quad \text{and}\\
            \tsum_{j=1}^n x_{i+n, j+n} c_j = \tsum_{j=1}^n \gamma_j c_j
        \end{align*}
        for $i \in \{1, \dots, n\}$.
\end{enumerate}
\end{proposition}

\begin{proof}
Let $\pi_\lambda : R \rightarrow R^+/\Ann_\lambda(R)$ and $\pi_\rho : R \rightarrow R^+/\Ann_\rho(R)$ be
the canonical homomorphisms. Then $\phi$ satisfies Condition~(1) of Theorem~\ref{max_ring_algebra_existence} if and
only if $\pi_\lambda(\phi(a_ia_j)) = \phi_1(\pi_\lambda(a_ia_j))$ and
$\pi_\rho(\phi(a_ia_j)) = \phi_2(\pi_\rho(a_ia_j))$. Since the action of~$\phi$ on $a_ia_j = f_R(b_i, c_j)$ is
defined by $f_R(\phi_1(b_i), c_j)$, we get
\begin{align*}
\pi_\lambda(f_R(\tsum_{\ell = 1}^n d_{i\ell}b_{\ell}, c_j)) &=
\pi_\lambda(\tsum_{\ell = 1}^n d_{i\ell} a_{\ell} a_j) =
\tsum_{\ell = 1}^n d_{i\ell} \tsum_{k=1}^n s_{\ell jk} b_{k} = \\
&= \tsum_{k=1}^n s_{ijk} \tsum_{\ell=1}^n d_{k \ell} b_\ell =
\tsum_{k=1}^n s_{ijk} \phi_1(b_k).
\end{align*}
Hence we have $\pi_\lambda(\phi(a_ia_j)) = \phi_1(\pi_\lambda(a_ia_j))$ if and only if $(d_{ij})$ is a solution of (i)
for $i,j \in \{1, \dots, n\}$. Similarly, $\pi_\lambda(\phi(a_ia_j)) = \phi_1(\pi_\lambda(a_ia_j))$ is satisfied
if and only if $(d_{ij})$ is a solution of (ii) for $i,j \in \{1, \dots, n\}$.

For every tuple $(d_{ij})$ we have $\sum_{j=1}^n d_{ij} b_j = \phi_1(\pi_\lambda(a_i))$ and
$\sum_{j=1}^n d_{i+n, j+n} c_j = \phi_2(\pi_\rho(a_i))$.
Let $\eta : R \rightarrow R^+/\Ann_\lambda(R) \times R^+/\Ann_\rho(R)$ be the canonical diagonal group homomorphism.
Then every tuple $(\pi_\lambda(r), \pi_\rho(r)) \in \eta(R)$ is of the form
$(\pi_\lambda(r), \pi_\rho(r)) = (\sum_{j=1}^n \gamma_j b_j, \sum_{j=1}^n \gamma_j c_j)$
with $\gamma_1, \dots, \gamma_n \in \mathbb{Z}$. We conclude that $(d_{i,j})$ is a
solution of the system in (2) if and only if
$(\phi_1(\pi_\lambda(a_i)), \phi_2(\pi_\rho(a_i))) = (\pi_\lambda(r), \pi_\rho(r))$ for some $r \in R$.
This is equivalent to $\phi$ satisfying Condition~(2) of Theorem~\ref{max_ring_algebra_existence}.
\end{proof}

Combining these results with those of the previous section, we get the following algorithm.

\begin{corollary}\label{Cor:maxring_algebra}
Let $R$ be a finite $\ZZ$-algebra given by a presentation as in Remark~\ref{remark:algebra_input}.
Then the following steps define an algorithm which returns elements $\phi_1, \dots, \phi_t$
of $\End(R^+/\Ann_\lambda(R)) \times \End(R^+/\Ann_\rho(R))$ corresponding to a system of $\ZZ$-module
generators of $\mathfrak{S}(R)$, together with an ideal $I \subseteq \ZZ[x_1, \dots, x_t]$ such that
we have a $\ZZ$-algebra presentation $\mathfrak{S}(R) = \ZZ[x_1, \dots, x_t]/I$.
\begin{enumerate}
\item[(1)] Using Remark~\ref{present_Ann}, compute presentations
of the $\ZZ$-modules $R^+/\Ann_\lambda(R)$ and $R^+/\Ann_\rho(R)$.

\item[(2)] For $i,j=1,\dots,n$, let $m_{ij}=\sum_{k=1}^n s_{ijk}a_k$.
Using Remark~\ref{solve_over_N}, solve the linear equation
$\sum_{i,j=1}^{n} x_{ij} m_{ij} =0$ in~$(\ZZ^n)^n$ and
deduce a finite presentation of the $\ZZ$-module $R^2$.

\item[(3)] Apply Proposition~\ref{mrs_endos} to compute tuples
$(\delta_{ij}^{(1)}), \dots,(\delta_{ij}^{(r)}) \in (\ZZ^{2n})^{2n}$
corresponding to a system of $\ZZ$-module generators of $\mathfrak{S}(f_R)$.

\item[(4)] Solve the system of homogeneous linear equations over
$R^+/\Ann_\lambda(R) \times R^+/\Ann_\rho(R)$ given by
\begin{align*}
    &\tsum_{\ell=1}^n \tsum_{k=1}^n s_{\ell j k} x_{i\ell}b_k
    - \tsum_{k=1}^n\tsum_{\ell=1}^n s_{ijk} x_{k\ell} b_{\ell} = 0 \quad \text{and} \\
    &\tsum_{\ell=1}^n \tsum_{k=1}^n s_{\ell j k} x_{i+n, \ell+n}c_k
    - \tsum_{k=1}^n\tsum_{\ell=1}^n s_{ijk} x_{k+n, \ell+n} c_{\ell} = 0
\end{align*}
for $i,j \in \{1, \dots, n\}$. Let $L_1 \subseteq (\ZZ^{2n})^{2n}$ be its solution space.

\item[(5)] Solve the system of homogeneous linear equations over
$R^+/\Ann_\lambda(R) \times R^+/\Ann_\rho(R)$ given by
\begin{align*}
    \tsum_{j=1}^n x_{ij} b_j - \tsum_{j=1}^n y_j b_j = 0 \quad \text{and}\\
    \tsum_{j=1}^n x_{i+n, j+n} c_j - \tsum_{j=1}^n y_j c_j = 0
\end{align*}
for $i \in \{1, \dots, n\}$. Let $L_2 \subseteq (\ZZ^{2n})^{2n}$ be the projection of
the solution space onto the $x$-coordinates.

\item[(6)] Compute the intersection of $L_1$, $L_2$ and the
$\ZZ$-module $\langle (\delta_{ij}^{(1)}), \dots,(\delta_{ij}^{(r)})\rangle$.
Let $(d_{ij}^{(1)}), \dots,(d_{ij}^{(t)})$ be tuples in $(\ZZ^{2n})^{2n}$ which
form a $\ZZ$-basis of this intersection.

\item[(7)] Apply Proposition~\ref{mrs_pres} to the tuples $(d_{ij}^{(1)}), \dots,(d_{ij}^{(t)})$
and obtain an ideal $I \subseteq \ZZ[x_1, \dots, x_t]$.

\item[(8)] Return the tuples of endomorphisms $\phi_1, \dots, \phi_t$ defined by
$\phi_k = (\phi_{k1}, \phi_{k2})$ with $\phi_{k1}(b_i) = \sum_{j=1}^n d_{ij}^{(k)} b_j$
and $\phi_{k1}(c_i) = \sum_{j=n+1}^{2n} d_{ij}^{(k)} c_{j-n}$, together with the ideal~$I$.

\end{enumerate}
\end{corollary}

\begin{proof}
An element $v \in (\ZZ^{2n})^{2n}$ of the intersection in Step~(6) corresponds to an
element~$\phi$ of $\End(R^+/\Ann_\lambda(R)) \times \End(R^+/\Ann_\rho(R))$ such that
$\phi \in \mathfrak{S}(f_R)$.
Since the element~$v$ is a solution of the system in Step~(4), it satisfies
Condition~(1) in Proposition~\ref{prop:linear_condition}. It is also a solution
of the system in Step~(5) and therefore satisfies Condition~(2) in
Proposition~\ref{prop:linear_condition}. Hence we get $\phi \in \mathfrak{S}(R)$.
This means that the intersection in Step~(6) represents exactly the elements of
$\mathfrak{S}(R)$. Therefore Step~(8) returns a $\ZZ$-algebra presentation of $\mathfrak{S}(R)$.
\end{proof}

Notice that all computations in this algorithm reduce to solving linear systems 
over a finitely presented abelian group.
Solutions of such linear systems can be obtained in polynomial time with respect 
to the size of the input, as described in
Remark~\ref{solve_over_N} and~\ref{linalgZcomplexity}.

\bigbreak
%
%

\section{Computing the Primitive Idempotents}%
\label{idem}

In the preceding sections we saw how to compute a $\ZZ$-algebra presentation
of the maximal ring of scalars $\mathfrak{S}(R)$ of a finite $\ZZ$-algebra $R$.
Now we want to decompose $\mathfrak{S}(R)$ into irreducible factors, i.e., we want to calculate
its primitive idempotents. Thus we use the following setting.

Let $S$ be a commutative ring which is a finitely generated $\ZZ$-module.
We assume that~$S$ is given by a presentation
$$
S \cong \ZZ[x_1,\dots,x_n]/I
$$
where~$I$ is an ideal in $P=\ZZ[x_1,\dots,x_n]$ which is given
by an explicit set of generators.

Recall that an idempotent~$e$ of~$S$ is called {\bf primitive} if it is non-zero
and not of the form $e=e'+e''$ with idempotents~$e',e''$ such that $e'e''=0$.
In order to compute the primitive idempotents of~$S$, we shall use
strong Gr\"obner bases which are defined as follows.

\begin{definition} Let $\mathbb{T}^n = \{x_1^{\alpha_1}\cdots x_n^{\alpha_n}
\mid \alpha_i \ge 0\}$  be the monoid of {\bf terms} in~$P$.
\begin{enumerate}
\item[(a)] A {\bf term ordering} $\sigma$ on~$\mathbb{T}^n$
is a complete ordering relation which is a well-ordering and compatible
with multiplication. In other words, for $t,t',t''\in\mathbb{T}^n$ we have
$1\le_\sigma t$ and $t\le_\sigma t'$ implies $t\,t'' \le_\sigma t'\, t''$.

\item[(b)] Every polynomial $f\in P \setminus \{0\}$ has a unique representation
$f=c_1t_1 + \cdots +c_st_s$ with $c_i\in \ZZ \setminus \{0\}$,
with $t_i\in\mathbb{T}^n$, and with $t_1 >_\sigma \cdots >_\sigma t_s$.
Then $\LM_\sigma(f)=c_1t_1$ is called the {\bf leading monomial} of~$f$,
the number $\LC_\sigma(f)=c_1$ is called its
{\bf leading coefficient}, and the term $\LT_\sigma(f)=t_1$ is called its
{\bf leading term} with respect to~$\sigma$.

\item[(c)] Given an ideal~$I$ in~$P$, a set of polynomials $G=\{g_1,\dots,g_r\}$
in~$I$ is called a {\bf strong $\sigma$-Gr\"obner basis} of~$I$ if, for every polynomial
$f\in I\setminus \{0\}$, there exists an index $i\in\{1,\dots,r\}$ such that $\LM_\sigma(f)$
is a multiple of $\LM_\sigma(g_i)$.
\end{enumerate}
\end{definition}

Strong Gr\"obner bases can be computed using a suitable extension of Buchberger's
Algorithm (see~\cite{AL}, Ch.~4). Many improvements for this computation have been found
and implemented (see for instance~\cite{Lic} and~\cite{Pop}). Subsequently,
our main task will be reduced to computing several strong Gr\"obner bases.
For some ideal-theoretic operations which can be performed effectively using
strong Gr\"obner bases, we refer to~\cite{AL}, Ch.~4 and~\cite{KR1}, Ch.~3.

The first step is to calculate the primary decomposition of the ideal~$I$.
For this, we use a simplification of the algorithm presented in~\cite{GTZ}
and~\cite{PSS}. Let us collect some basic observations about the structure
of the primary decomposition of~$I$.

\begin{remark}\label{primary_props}
Let $I$ be an ideal in $P=\ZZ[x_1,\dots,x_n]$
such that $P/I$ is a finitely generated $\ZZ$-module.
\begin{enumerate}
\item[(a)] The ring $S=P/I$ satisfies $\dim(S) \le 1$. Hence the 
primary components of~$I$ have height~$n$ or~$n+1$.

\item[(b)] A primary component~$\mathfrak{Q}$ of~$I$ has height~$n+1$ if and only if
it contains an integer $q\ge 2$. In this case, the radical of~$\mathfrak{Q}$
is a maximal ideal~$\mathfrak{M}$ of~$P$. 
Then we have $\mathfrak{M}\cap \ZZ = \langle p\rangle$
for some prime number~$p$ and $q=p^e$ for some $e>0$.
The residue class ideal~$\mathfrak{m}$ of~$\mathfrak{M}$ in
$\FF_p[x_1,\dots,x_n]$ is a maximal ideal, too.

\item[(c)] In the setting of~(b), we may recover the ideal~$\mathfrak{M}$
from~$\mathfrak{m}$ as follows. For $\bar{a}\in\FF_p$, we call the unique
element $a\in\{0,\dots,p-1\}$ such that $\bar{a}=a+\langle p\rangle$ the
{\it canonical lifting} of~$\bar{a}$. For every polynomial in $\FF_p[x_1,\dots,x_n]$,
we define its {\it canonical lifting} in~$P$ by lifting  all coefficients canonically.
Given a system of generators $\{\bar{g}_1,\dots,\bar{g}_s\}$ of~$\mathfrak{m}$,
we then have $\mathfrak{M} = \langle p, g_1,\dots,g_s\rangle$, where~$g_i$ is
the canonical lifting of~$\bar{g}_i$ for $i=1,\dots,s$.

\item[(d)] A primary component~$\mathfrak{Q}$ of~$I$ has height~$n$
if and only if $\mathfrak{Q} \cap \ZZ = \{0\}$. In this case
we have $\mathfrak{Q}= \mathfrak{Q} \QQ[x_1,\dots,x_n]\cap P$, 
and the extended ideal $\mathfrak{Q}\QQ[x_1,\dots,x_n]$ is primary to
a maximal ideal.
\end{enumerate}
\end{remark}

The next proposition allows us to split the task of computing the primary
composition of~$I$ into two cases. Recall that a strong Gr\"obner basis is
called {\bf minimal} if no leading monomial divides another one.

\begin{proposition}\label{lemma_Q}
Let $I$ be an ideal in $P=\ZZ[x_1,\dots,x_n]$ such that
$P/I$ is a finite $\ZZ$-algebra, let $G=\{g_1,\dots,g_s\}$
be a minimal strong Gr\"obner basis of~$I$, and let $N\in\ZZ$
be the least common multiple of the leading coefficients of the
elements of~$G$.
\begin{enumerate}
\item[(a)] The ideal $(I : \langle N\rangle)/I$ is the torsion
subgroup of~$P/I$.

\item[(b)] We have $I = (I : \langle N\rangle) \cap (I+ \langle N\rangle)$.
\end{enumerate}
\end{proposition}

\begin{proof}
To prove~(a), it obviously suffices to show the inclusion $\supseteq$.
Let $f\in P$ be a polynomial such that $f+I$ is in the torsion of~$P/I$.
By~\cite{AL}, Prop.~4.4.4 and Ex.~4.4.1, we know that $N^i f\in I$ for some
$i\ge 1$. We want to show that $Nf\in I$. Assuming that this is not the
case, there exists a smallest number $i\ge 1$ such that $N^i f\notin I$ and
$N^{i+1}f\in I$. Among all polynomials~$f$ with this property, choose the
one with the smallest leading term w.r.t.~$\sigma$. Then $N^{i+1}f\in I$
implies that there exists a $j\in\{1,\dots,s\}$ such that $\LM_\sigma(g_j)$
divides $\LM_\sigma(N^{i+1}f)$. Write $\LM_\sigma(g_j)=c_jt_j$ with $c_j\in\ZZ$
and $t_j\in\mathbb{T}^n$. Then~$t_j$ divides $\LT_\sigma(f)$ and we can write
$\LT_\sigma(f)=t_j\,\tilde{t}$ with $\tilde{t}\in\mathbb{T}^n$. Since~$c_j$
divides~$N$, it follows that  $\LM_\sigma(g_j)$ divides $\LM_\sigma(Nf)$.
Hence there is a monomial~$h$ such that $Nf-hg_j$ has a smaller leading term
than~$f$. Now $Nf-hg_j\notin I$ and $N^i(NF-hg_j)\in I$ contradict the
minimality of~$\LT_\sigma(f)$.

Since~(a) implies $I:\langle N\rangle = I: \langle N\rangle^\infty$,
claim~(b) is a standard lemma in commutative algebra. In fact, if
$f\in P$ with $Nf\in I$ and $f=g+Nh$ for some $g\in I$, $h\in P$
then $N^2h\in I$ implies $Nh\in I$, and hence $f\in I$.
\end{proof}

Note that the number $N$ from the previous proposition is in general not the exponent of $P/I$.

\begin{example}
Consider the ideal $I = \langle x^2,y^2,z^2,xz+yz,xy,2x-y, 3z \rangle \subseteq \ZZ[x,y,z]$.
The generators form a strong Gr\"obner basis of $I$ and $P/I$ is generated as a $\ZZ$-module
by the residue classes of $x,y,z$. Even though $2x \in \LM(I)$, the element $\bar{x}$ is not
in the torsion of $P/I$, since the normal form of $2x$ with respect to~$I$ is~$y$, 
and no proper multiple of~$y$ is in~$I$. Therefore the exponent of $P/I$ is~3.
\end{example}

At this point we are ready to formulate Algorithm~\ref{decomp}
for computing the primary decomposition of~$I$.

\begin{algorithm} \caption{Computing a Primary Decomposition}\label{decomp}
\begin{algorithmic}[1]
\REQUIRE Generators of an ideal $I \subseteq P$ such that $P/I$ is a finite $\ZZ$-algebra.
\ENSURE A tuple of ideals~$L=(\mathfrak{Q}_1, \dots, \mathfrak{Q}_k)$ such that
        $I = \mathfrak{Q}_1 \cap \dots \cap \mathfrak{Q}_k$ is a primary decomposition of~$I$.
   \STATE Compute a minimal strong Gr\"obner basis $\{g_1, \dots, g_s\}$ of~$I$.
   \STATE Let $q$ be the generator of $I \cap \ZZ$.
   \STATE $L := (\;)$
   \IF {$q = 0$}
       \STATE $N := \lcm(\LC(g_1), \dots, \LC(g_s))$
       \STATE Compute a primary decomposition $I\,\QQ[x_1, \dots, x_n]=
               \overline{Q}_1 \cap \dots \cap \overline{Q}_\ell$.
       \STATE Compute $\overline{Q}_j \cap P$ and append these ideals to~$L$.
       \STATE Recursively apply the algorithm to $I + \langle N \rangle$
              and obtain the set~$M$.
       \STATE Compute $J := \bigcap_{{\scriptstyle\mathfrak{Q}}\, \in L} \mathfrak{Q}$.
       \STATE Remove all ideals in~$M$ that contain~$J$.
       \RETURN $L \cup M$
    \ELSE
       \STATE Compute the prime factorization $q=p_{1}^{\nu_1}\cdots p_{r}^{\nu_r}$.
       \STATE $M := (\;)$
       \FOR{$i=1,\dots,r$}
          \STATE Compute a primary decomposition $I\,\FF_{p_i}[x_1, \dots, x_n]
                 = \overline{Q}_1 \cap \dots \cap \overline{Q}_m$.
          \IF{$\nu_i=1$}
             \STATE Compute the preimages~$\mathfrak{Q}_j$ of~$\overline{Q}_j$ in~$P$
                    and append them to~$M$.
          \ELSE
             \STATE Compute the prime components $\overline{P}_1, \dots,
                    \overline{P}_m$\/ of $I\,\FF_{p_i}[x_1, \dots, x_n]$.
             \STATE For $j=1,\dots,m$, compute the preimage~$\mathfrak{P}_j$ 
                    of~$\overline{P}_j$ in~$P$.
             \STATE For $j=1,\dots,m$, calculate $k \in \mathbb{N}$ such that
             $(I+\mathfrak{P}^k) \cap (I:\mathfrak{P}^\infty) = I$.
             \STATE For $j=1,\dots,m$, append the ideal $\mathfrak{Q}_j := I+\mathfrak{P}^k$ to~$M$.
         \ENDIF
      \ENDFOR
      \RETURN $M$
   \ENDIF
\end{algorithmic}
\end{algorithm}

Let us point out how the various steps of this algorithm
can be carried out effectively.

\begin{remark} Let $I$ be an ideal in~$P$ such that $P/I$ is a finitely
generated $\ZZ$-module.
\begin{enumerate}
\item[(a)] The ideal $I\,\QQ[x_1, \dots, x_n]$ in Line~6 is a 0-dimen\-sional ideal
in $\QQ[x_1,\dots,x_n]$. Hence its primary decomposition can be computed using
standard computer algebra methods (see for instance~\cite{KR3}, Ch.~5).

\item[(b)] If $I$ contains a non-zero integer~$q$ and~$p$ is a prime dividing~$q$,
then also $I\,\FF_p[x_1, \dots, x_n]$
is a 0-dimen\-sional ideal. We can therefore again use standard primary decomposition
algorithms for 0-dimensional ideals in Line~16 such as the ones in~\cite{KR3}, Ch.~5.

\item[(c)] Notice that most of the algorithms mentioned in~(a) and~(b) also
determine the {\bf prime components} of the respective ideals, i.e., the radical
ideals of the primary components.

\item[(d)] The preimages in Line~21 can be found as in Remark~\ref{primary_props}.b.

\end{enumerate}
\end{remark}

The following results, adapted from~\cite{PSS}, are used to prove the correctness of
Algorithm~\ref{decomp}. They allow us to reduce the task to the computation of primary
decompositions in $\QQ[x_1, \dots, x_n]$ and $\FF_p[x_1, \dots, x_n]$.

\begin{lemma}\label{lemma_liftings}
Let $I \subseteq P$ be an ideal such that $S=P/I$ is a finite $\ZZ$-algebra.
\begin{enumerate}
\item[(a)] Suppose that $I \cap \ZZ = \{0\}$, and
let $I\,\QQ[x_1, \dots, x_n]= \overline{Q}_1 \cap \dots \cap \overline{Q}_s$
be a primary decomposition. If we set $\mathfrak{Q}_i = \overline{Q}_i \cap P$ 
for $i=1, \dots, s$, then $I = \mathfrak{Q}_1 \cap \cdots \cap \mathfrak{Q}_s$ 
is a primary decomposition of~$I$.

\item[(b)] Suppose that $I \cap \ZZ = \langle p \rangle$ for some prime number~$p$,
and let $I\,\FF_p[x_1, \dots, x_n]= \overline{Q}_1 \cap\cdots\cap \overline{Q}_s$
be a primary decomposition.  If we let $\mathfrak{Q}_i$ be the preimage 
of~$\overline{Q}_i$ in~$P$ for $i=1,\dots,s$, then
$I=\mathfrak{Q}_1 \cap \cdots \cap \mathfrak{Q}_s$ is a primary decomposition of~$I$.

\item[(c)] Suppose that $I \cap \ZZ = \langle p^\nu \rangle$ for some prime number~$p$ 
and $\nu\ge 1$, and let $\overline{P}_1, \dots, \overline{P}_s$ 
be the prime components of $I\,\FF_p[x_1, \dots, x_n]$. If we let $\mathfrak{P}_i$ be the
preimage of~$\overline{P}_i$ in~$P$ for $i=1,\dots,s$, then
$\mathfrak{P}_1, \dots, \mathfrak{P}_s$ are the prime components of~$I$.
\end{enumerate}
\end{lemma}

\begin{proof}
First we prove~(a). Consider the ring homomorphism $\Phi:\; \ZZ[x_1, \dots, x_n]
\longrightarrow \QQ[x_1, \dots, x_n]$ given by embedding the coefficients
into~$\QQ$. 

Then we have
\begin{align*}
 I \;&=\;  I\,\QQ[x_1,\dots,x_n] \cap P \;=\; \Phi^{-1}(I\,\QQ[x_1, \dots, x_n]) \\
   \;&=\; \Phi^{-1}(\overline{Q}_1) \cap \cdots \cap \Phi^{-1}(\overline{Q}_s) 
      \;=\; \mathfrak{Q}_1 \cap \cdots \cap \mathfrak{Q}_s.
\end{align*}
Since the preimage of a primary ideal is primary, this implies the claim.

To prove~(b), we proceed analogously. We consider the ring
homomorphism $\Psi:\; \ZZ[x_1, \dots, x_n] \longrightarrow
\FF_p[x_1, \dots, x_n]$ given by $\Psi(c_1 t_1 + \cdots + c_m t_m) =
\bar{c}_1 t_1 + \cdots + \bar{c}_m t_m$ for $c_i \in \ZZ \setminus \{0\}$
and $t_i \in \mathbb{T}^n$ and argue as before.

Finally, we show~(c). If $p^\nu \in I$ then $p \in \sqrt{I}$, and we can apply~(b)
to compute the primary components of~$\sqrt{I}$. Since~$\sqrt{I}$ is a radical ideal,
these are prime ideals and thus precisely the prime components of~$I$.
\end{proof}

Now we can prove the correctness of Algorithm~\ref{decomp}.

\begin{proposition}\label{prop:decomp}
Algorithm~\ref{decomp} terminates and computes a primary decomposition
of the input ideal~$I$.
\end{proposition}

\begin{proof}
Since the first recursive call of the algorithm necessarily leads to the else branch,
the algorithm terminates. It remains to prove the second claim.

By Proposition~\ref{lemma_Q}, we have $I: \langle N\rangle = I : \langle N\rangle^\infty$
and $I= (I: \langle N\rangle) \cap (I +\langle N\rangle)$ for the number
$N>0$ computed in Step~(5). Hence, in the case $I \cap\ZZ=\{0\}$, we may
compute the primary decomposition of
$$
I\;:\; \langle N\rangle \;=\;  I\,\QQ[x_1, \dots, x_n] \cap P
$$
via Lemma~\ref{lemma_liftings}.a
and the ideal $I + \langle N \rangle$ can be decomposed by the following steps,
since it intersects~$\ZZ$ non-trivially.

If $I \cap \ZZ = \langle q \rangle$ for some $q>0$,
we determine the prime factorization $q=p_1^{\nu_1}\cdots p_r^{\nu_r}$
in Step~(13). Then we have $I=\bigcap_{i=1}^r (I + \langle p_i^{\nu_i} \rangle )$,
and thus it suffices to decompose each ideal $ I + \langle p_i^{\nu_i} \rangle$
individually.

In the case $\nu_i=1$, we may again use Lemma~\ref{lemma_liftings}
and conclude that it is enough to compute the primary decomposition of
$I\,\FF_{p_i}[x_1, \dots, x_n]$ and to calculate the preimages
of the primary components in~$P$. In the case $\nu_i > 1$,
we can still compute a set of minimal associated primes of~$I$.
They are the preimages of the minimal associated primes of $I\,\FF_{p_i}[x_1,\dots,x_n]$
in~$P$. Then we can extract the primary components $\mathfrak{Q}_j$ from the corresponding
prime components~$\mathfrak{P}_j$. Since $\mathfrak{P}_j$ is a maximal ideal,
$I+\mathfrak{P}_j^k$ is $\mathfrak{P}_j$-primary for every positive $k$. If $k$ is sufficiently
large, $I+\mathfrak{P}_j^k$ is a primary component of $I$ (see \cite{EHV}, Section~4). We can
determine such a $k$ by checking whether $(I+\mathfrak{P}_j^k) \cap (I: \mathfrak{P}_j) = I$
(see \cite{IsYo}, Criterion 3).

Altogether, the collection of all primary components in~$L \cup M$ intersects to~$I$.
If $\bigcap_{i \neq j} \mathfrak{Q}_j \subset \mathfrak{Q}_i$ for some primary ideal
$\mathfrak{Q}_i \in L \cup M$ then~$\mathfrak{Q}_i$ has to be an element of~$M$.
But this means the condition $J \subseteq Q_i$ in Step~(10) is satisfied.
Hence we have found a primary decomposition of~$I$.
\end{proof}

A small modification of Algorithm~\ref{decomp} allows us to compute the prime components
of the ideal~$I$ as well.

\begin{remark}\label{prim_comp}
In Algorithm~\ref{decomp}, replace Steps~(6) and (20) by
\begin{enumerate} 
\item[(6')] Compute a primary decomposition $I\, \QQ[x_1,\dots,x_n] = \overline{Q}_1
\cap \cdots \cap \overline{Q}_\ell$ and the corresponding prime components
$\overline{P}_1,\dots,\overline{P}_\ell$.

\item[(20')] Compute a primary decomposition $I\, F_{p_i}[x_1,\dots,x_n] =
\overline{Q}_1 \cap \cdots \cap \overline{Q}_m$ and the corresponding prime components
$\overline{P}_1,\dots,\overline{P}_m$.
\end{enumerate}
Furthermore, include the computation of $\mathfrak{P}_i=\overline{Q}_i \cap P$ in Step~(7),
include the computation of the preimages $\mathfrak{P}_i$ of~$\overline{P}_i$ in Steps~(18)
and~(21), and append the pairs $(\mathfrak{Q}_i,\mathfrak{P}_i)$ to the lists~$L$
and~$M$ in the appropriate places. 

Then the resulting algorithm computes a primary
decomposition $I=\mathfrak{Q}_1\cap \cdots \cap \mathfrak{Q }_k$ of~$I$ as well as
the corresponding prime components~$\mathfrak{P}_i$. In particular, we have a primary
decomposition $\sqrt{I}= \mathfrak{P}_1 \cap\cdots\cap \mathfrak{P}_k$.
\end{remark}

Having calculated a primary decomposition $I=\mathfrak{Q}_1\cap \cdots \cap \mathfrak{Q}_k$ 
of the ideal~$I$, our next goal is to find the primitive idempotents of $S=P/I$.
Recall that the {\bf spectrum} of the ring~$S$ is the set $\Spec(S)$ of all prime ideals
of~$S$. This set is a topological space via the {\bf Zariski topology} whose closed
sets are the sets $\mathcal{V}(J)$ consisting of all prime ideals containing some ideal~$J$ of~$S$.
The spectrum and the primitive idempotents of~$S$ are related through the following observations
(see~\cite{Bou}, II.4.3, ex. 14 and~\cite{PS}, ex. 4.B.5).

\begin{remark}
Let~$T$ be a commutative, unitary, noetherian ring. 
\begin{enumerate}
\item[(a)] Given an idempotent $e\in T$, the set $\mathcal{V}(1-e)$
is both open and closed in~$\Spec(T)$.

\item[(b)] If $U \subseteq \Spec(T)$ is a subset which is both open and closed,
there exists a unique idempotent $e\in T$ such that in $T_{\mathfrak{p}}/\mathfrak{p}
T_{\mathfrak{p}}$ we have $\bar{e}=1$ for $\mathfrak{p}\in U$ and $\bar{e}=0$ otherwise.

\item[(c)] The correspondence given in~(a) and~(b) is 1-1. The primitive idempotents
correspond uniquely to the connected components of~$\Spec(T)$.
\end{enumerate}
\end{remark}

Thus, in order to compute the primitive idempotents of~$S$, we calculate the connected
components of~$\Spec(S)$ first. If~$\dim(S)=1$, this ring has infinitely many prime ideals.
To describe the connected components of~$\Spec(S)$, the following definition will come in handy.

\begin{definition}
Let $I \subseteq P$ be an ideal, and let $I = \mathfrak{Q}_1 \cap \cdots \cap \mathfrak{Q}_k$ 
be a primary decomposition of~$I$. A maximal subset of $\{\mathfrak{Q}_1, \dots, \mathfrak{Q}_k\}$ 
such that all corresponding prime ideals are part of the same connected component 
of $\Spec(P/I)$ is called a \textbf{connected component} of the primary decomposition.
\end{definition}

Recall that, given a primary decomposition $I=\mathfrak{Q}_1 \cap \cdots\cap \mathfrak{Q}_k$
of an ideal $I\subseteq P$ such that $S=P/I$ is a finite $\ZZ$-algebra,
the primary components $Q_i$ are either of height~$n$ and do not contain
a non-zero integer, or they are of height $n+1$ and correspond to a maximal component of~$I$.
The connected components of the primary decomposition of~$I$ are determined by
the next Algorithm~\ref{alg_connected}.

\begin{algorithm}\caption{Computing the Connected Components}\label{alg_connected}
\begin{algorithmic}[1]
\REQUIRE An ideal $I \subseteq P$ such that $P/I$ is a finite $\ZZ$-algebra.
\ENSURE A set $M=\{C_1,\dots,C_\nu\}$ such that $C_1,\dots,C_\nu$ are the connected components 
        of a primary decomposition of~$I$.
  \STATE Using Remark~\ref{prim_comp}, compute a primary decomposition of~$I$ together
         with the corresponding prime components. Let $\mathfrak{Q}_1, \dots, \mathfrak{Q}_\ell$ 
         be the primary components of height~$n+1$, let $\mathfrak{M}_1,\dots,\mathfrak{M}_\ell$
         be the corresponding prime components, let $\mathfrak{Q}'_1, \dots, \mathfrak{Q}'_m$ 
         be the primary components of height~$n$, and let $\mathfrak{P}_1,\dots,\mathfrak{P}_m$
         be the corresponding prime components.
    \STATE Let $M = \{ \{\mathfrak{Q}'_1\}, \dots, \{\mathfrak{Q}'_m\}\}$.
    \FOR{$i=1$ \TO $m$}
        \FOR{$j=1$ \TO $\ell$}
        \STATE If $\mathfrak{P}_i \subseteq \mathfrak{M}_j$ then append $\mathfrak{Q}_j$ 
             to the set in~$M$ that contains $\mathfrak{Q}'_i$.
        \ENDFOR
    \ENDFOR
    \WHILE{there are sets $C,C' \in M$ such that there exist
           $\mathfrak{Q}'_i \in C$ and $ \mathfrak{Q}'_j \in C'$ with $\mathfrak{P}_i + \mathfrak{P}_j 
           \ne \langle 1 \rangle$}
    \STATE replace~$C$ and~$C'$ in~$M$ by $C \cup C'$.
    \ENDWHILE
    \STATE For every ideal $\mathfrak{Q}_i$ which is not contained in any of the sets of~$M$, append
           the set $\{ \mathfrak{Q}_i \}$ to~$M$.
    \RETURN $M$.
\end{algorithmic}
\end{algorithm}

\begin{proposition}
Algorithm~\ref{alg_connected} computes the connected components of a primary decomposition of~$I$.
\end{proposition}

\begin{proof}
In view of Remark~\ref{primary_props}, a connected component~$C$ of the primary decomposition
$I=\mathfrak{Q}_1\cap \cdots \cap \mathfrak{Q}_k \cap \mathfrak{Q}'_1 \cap \cdots \cap 
\mathfrak{Q}'_m$ computed in Step~(1) is of one of the following forms.
\begin{enumerate}
\item[(a)] There exists a subset~$T$ of $\{\mathfrak{Q'}_1,\dots,\mathfrak{Q}'_m\}$
such that~$C$ consists of~$T$ and all ideals~$\mathfrak{Q}_j$ containing an ideal in~$T$.

\item[(b)] We have $C = \{ \mathfrak{Q}_i \}$ for some $i\in \{1,\dots,\ell\}$.
\end{enumerate}
Here two ideals $\mathfrak{Q}'_i,\mathfrak{Q}'_j$ are in the same connected component
if and only if there exists a maximal ideal containing both~$\mathfrak{P}_i$ and~$\mathfrak{P}_j$,
and this is equivalent to $\mathfrak{P}_i + \mathfrak{P}_j \ne \langle 1 \rangle$.
Thus the ideals $\mathfrak{Q}'_i$ in connected components of type~(a) are correctly determined 
in Steps~(8)-(10), and the ideals~$\mathfrak{Q}_i$ in these components are sorted out
in the loop (3)-(7). 

The remaining ideals $\mathfrak{Q}_i$ correspond to connected components
$\{\mathfrak{M}_i\}$ of~$\Spec(P/I)$ for which~$\mathfrak{M}_i$ is both a minimal and a 
maximal ideal, i.e., to isolated points of~$\Spec(P/I)$. They are correctly found in Step~(11).
\end{proof}

Finally, we can calculate the desired idempotents of~$S$ as follows.

\begin{proposition}{\bf (Computing the Primitive Idempotents)}\label{alg_idempotents}\\
Let $I$ be an ideal in~$P$ such that $S = P/I$ is a finite $\ZZ$-algebra.
The following steps define an algorithm which computes the primitive idempotents of~$S$.
\begin{enumerate}
\item[(1)] Using Algorithm~\ref{alg_connected}, compute the connected components $C_1, \dots, C_\nu$ 
of a primary decomposition of~$I$.

\item[(2)] For $i=1,\dots,\nu$, compute $J_i = \bigcap_{Q \in C_i} Q$.

\item[(3)] For $i=1,\dots,\nu$, compute elements $q_i \in \bigcap_{j \ne i} J_j$
and $p_i \in J_i$ such that $q_i+p_i = 1$.

\item[(4)] Return $\{q_1, \dots, q_\nu\}$.
\end{enumerate}
\end{proposition}

\begin{proof}
First we show that the ideals $J_1, \dots, J_\nu$ are pairwise comaximal. It suffices to show
that all ideals $\mathfrak{Q} \in C_i$ and $\mathfrak{Q}'\in C_j$ are pairwise comaximal for $i \ne j$.
Since~$\mathfrak{Q}$ and $\mathfrak{Q}'$ are primary ideals, we may replace them by their 
prime components $\mathfrak{P}$ and~$\mathfrak{P}'$. Then we have 
$\mathfrak{P} + \mathfrak{P}' =  \langle 1 \rangle$, since
$\mathfrak{P}$ and~$\mathfrak{P}'$ would otherwise belong to the same connected 
component of~$\Spec(P/I)$.

Next we consider the canonical ring homomorphism
$$
\Phi \colon \; P/I \;\longrightarrow\; \prod_{i=1}^\nu P/J_i
$$
given by $\Phi(f+I)=(f+J_1,\dots,f+J_k)$ for $f\in P$. Since the ideals~$J_i$ are
pairwise comaximal, and since $J_1\cap \cdots \cap J_\nu =I$, the Chinese Remainder
Theorem (cf.~\cite{KR1}, Lemma~3.7.4) shows that~$\Phi$ is an isomorphism.
Furthermore, by the definition of the ideals~$J_i$, the spectrum $\Spec(P/J_i)$ is connected,
and therefore the ring $P/J_i$ has only the trivial idempotents.
The primitive idempotents of $P/I$ are then given by the preimages of the elements
$e_1,\dots,e_\nu$ on the right-hand side. The polynomials $q_i$ computed in Step~(3)
are precisely the preimages of the elements~$e_i$ for $i=1,\dots,\mu$.
\end{proof}

Let us finish with some observations on the complexity of the algorithms in this section.

\begin{remark}
Given a presentation of a commutative finite $\ZZ$-algebra $S=P/I$ as at the
beginning of this section, the computation of the primitive idempotents
of~$S$ is dominated by the calculation of a strong Gr\"obner basis of~$I$
in Line~1 of Algorithm~\ref{decomp}. Although a precise complexity estimate
seems not to have been derived yet, it appears to be singly exponential
in the number and size of the generators of~$S$.

The Gr\"obner basis computation can be avoided if a presentation of $S$ is given
as in Proposition~\ref{mrs_pres}. This means we are given $\ZZ$-module generators
and structure constants, together with a presentation
$S \cong \ZZ^m / U$ of~$S$ as a $\ZZ$-module, where $U$ is a submodule of $\ZZ^m$.
In a forthcoming paper~\cite{KW} we give a detailed complexity analysis of the computation 
of the primitive idempotents, if the $\ZZ$-algebra is given as above. 
We show that all calculations, except for one integer prime factorization in
Step~13 of Algorithm~\ref{decomp}, can be performed in probabilistic polynomial time 
in the bit complexity of the input.
\end{remark}

\bigbreak
%
%

\section{The Canonical Decomposition of a Bilinear Map}
\label{sec:decomp_bilinear}

Given a finite $\ZZ$-algebra $R$, our main goal is to use the decomposition
of its maximal ring of scalars~$\mathfrak{S}(R)$ into directly indecomposable
factors to get a decomposition of $R$ itself.
In this section we start by studying the decomposition of bilinear maps. Later we
apply the results to the multiplication maps of finite $\ZZ$-algebras.

\begin{definition}
Let $N_1, N_2 ,M$ be abelian groups, and let $f:\; N_1 \times N_2 \longrightarrow M$ be a bilinear map.
\begin{enumerate}
\item[(a)] If there exist subgroups~$M_i$ of~$M$, ~$N_{1i}$ of~$N_1$ and $N_{2i}$ of $N_2$ such that
    $M=\bigoplus_{i=1}^k M_i$, $N_1=\bigoplus_{i=1}^k N_{1i}$ and $N_2=\bigoplus_{i=1}^k N_{2i}$, and if
    there are bilinear maps $f_i:\;  N_{1i} \times N_{2i} \longrightarrow M_i$ for $i=1,\dots,k$ such that
$f(x_1 + \dots + x_k ,y_1 + \dots + y_k) = f_1(x_1, y_1) + \cdots + f_k(x_k,y_k)$
for all $x_i \in N_{1i}$ and $y_i \in N_{2i}$, then we say that~$f$ is the {\bf direct product}
of the maps~$f_i$, and we write $f = f_1 \times \cdots \times f_k$.

\item[(b)] The map $f$ is called \textbf{directly indecomposable} if~$f$ cannot
be written as a direct product of non-trivial bilinear maps.
\end{enumerate}
\end{definition}

A direct product decomposition of a bilinear map~$f$ into directly indecomposable maps
can be computed using the primitive idempotents of its maximal ring of scalars.

\begin{proposition}\label{decomp_bilinear}
Let $f:\; N_1 \times N_2 \longrightarrow M$ be a full, non-degenerate bilinear map,
and assume that we are given a complete set of primitive idempotents $\{e_1, \dots, e_k\}$
of the maximal ring of scalars~$\mathfrak{S}(f)$. Then we obtain direct sum decompositions
$N_1 = \bigoplus_{i=1}^{k} e_i N_1$, $N_2 = \bigoplus_{i=1}^{k} e_i N_2$ and $M = \bigoplus_{i=1}^{k} e_i M$
such that~$f$ is the direct product of the full, non-degenerate and directly indecomposable bilinear maps
$$
f_i:\;  e_i N_1 \times e_i N_2  \longrightarrow e_i M
$$
given by $f_i(e_i x,e_i y) = e_i f (x,y)$ for $x\in N_1$, $y\in N_2$ and for $i=1,\dots,k$.
\end{proposition}

\begin{proof}
Writing $1=e_1+\cdots + e_k$ in~$\mathfrak{S}(f)$, we see that
$x=1\cdot x = e_1x + \cdots +e_kx$ for every $x$ in $N_1$, $N_2$ or $M$. Thus we obtain the sum decompositions
$N_1=\sum_{i=1}^k e_i N_1$, $N_2=\sum_{i=1}^k e_i N_2$ as well as $M=\sum_{i=1}^k e_i M$. These sum decompositions
are direct, since for instance $e_i m = e_j m'$ with $m,m'\in M$ and $i\ne j$ implies
$e_i m = e_i^2 m = e_i e_j m' = 0$.

Clearly, for $i=1,\dots,k$, the maps $f_i:\; e_i N_1 \times e_i N_2 \longrightarrow e_i M$
defined by letting $f_i(e_i x,e_i y) = e_i f (x,y)$ for $x\in N_1$ and $y\in N_2$ are well-defined and bilinear.
To show that the map~$f_i$ is full, let $x \in e_iM$. Then there
exists an element $m \in M$ with $x = e_i m$. Since the map~$f$ is full, we find
elements $n_1, \dots, n_s \in N_1$ and $n'_1, \dots, n'_s \in N_2$
such that $m = f(n_1, n'_1) + \cdots + f(n_s, n'_s)$. Therefore we obtain
$$
x \;=\; e_i m \;=\; e_i^2 f(n_1, n'_1) + \cdots + e_i^2 f(n_s, n'_s)
  \;=\; f_i(e_i n_1, e_i n'_1) + \cdots + f_i(e_i n_s, e_i n'_s).
$$
To show the non-degeneracy of~$f_i$, assume that $x\in N_1$ satisfies $f_i(e_i x, e_i y)
= e_i f(x, y) = 0$ for all $y \in N_2$. Then we have $f(x,y) = 0$ for all $y \in N_2$,
since~$M$ is a faithful $\mathfrak{S}(f)$-module. Now the fact that~$f$ is non-degenerate implies $x=0$.

Finally, given elements $x\in N_1$ and $y\in N_2$, we let $x_i=e_i x$ and $y_i=e_iy$ for $i=1,\dots,k$,
and we calculate
$$
f(x,y) \;=\; f(x_1 + \cdots + x_k,\, y_1 + \cdots + y_k) \;=\; f_1(x_1, y_1) + \cdots + f_k(x_k, y_k)
$$
since $f(e_i x, e_j y) = f(x, e_i e_j y) = f(x, 0) = 0$ for all $i \ne j$. Thus we obtain the
direct product decomposition $f=f_1 \times \dots \times f_k$.

It remains to show that the maps~$f_i$ are directly indecomposable.
Suppose there exists a non-trivial decomposition $f_i = g \times h$.
Then we have $\mathfrak{S}(f_i) = \mathfrak{S}(g) \times \mathfrak{S}(h)$
by~\cite{Mya1}, Proposition~3.1, and we
get $e_i = e_g + e_h$ with orthogonal idempotents $e_g = (1,0)$ and $e_h = (0,1)$. This
contradicts the assumption that~$e_i$ is primitive.
\end{proof}

\begin{definition}
The decomposition $f= f_1 \times \cdots \times f_k$ given in this proposition
is called the {\bf canonical decomposition} of the bilinear map~$f$.
Using the results of the preceding sections, it can be computed as follows.
\end{definition}

\begin{corollary}{\bf (Computing the Canonical Decomposition of a Bilinear Map)}\label{cor:bilinear_decomp}
Let $f \colon N_1 \times N_2 \longrightarrow M$ be a full and non-degenerate bilinear map.
Assume that $\{a_1, \dots, a_n\}$, $\{a_{n+1}, \dots, a_{n'}\}$ and $\{b_1, \dots, b_m\}$ are generating sets of~$N_1$,
$N_2$ and~$M$, respectively. Furthermore, assume that we are given presentations of~$N_1$, $N_2$ and~$M$ as in
Remark~\ref{remark:alg_input} and structure constants $s_{ijk}\in \ZZ$ for $i=1,\dots,n$, $j=n+1, \dots, n'$ and
$k=1, \dots, m$ such that $f(a_i, a_j) = \sum_{k=1}^m s_{ijk} b_k$.
Then the following steps define an algorithm which computes the canonical
decomposition $f= f_1 \times \cdots \times f_\ell$.
\begin{enumerate}
\item[(1)] Use Proposition~\ref{mrs_endos} to compute tuples $(d_{ij}^{(1)}),\dots,(d_{ij}^{(r)}) \in (\ZZ^{n'})^{n'}$
corresponding to elements of~$\End(N_1) \times \End(N_2)$ which generate the maximal ring of
scalars~$\mathfrak{S}(f)$. For $k=1, \dots, r$ let $A_k \in \Mat_{n,n}(\ZZ)$ be the matrix given by
$(d_{i,j}^{(k)})_{i,j=1, \dots, n}$ and let ${A'}_k \in \Mat_{n'-n, n'-n}(\ZZ)$ be the matrix
given by $(d_{i,j}^{(k)})_{i,j=n+1, \dots, n'}$.

\item[(2)] Use Proposition~\ref{mrs_pres} to compute generators $g_1,\dots,g_k \in \ZZ[x_1,\dots,x_r]$
of an ideal~$I$ such that $\mathfrak{S}(f) = \ZZ[x_1,\dots, x_r]/I$.

\item[(3)] Use Proposition~\ref{alg_idempotents} to compute polynomials representing
a complete set $\{e_1, \dots, e_\ell\}$ of primitive idempotents of $\ZZ[x_1,\dots,x_r]/I$.

\item[(4)] For $i=1,\dots,\ell$, let $E_i = e_i(A_1, \dots, A_r)$ and $E'_i = e_i(A'_1, \dots, A'_r)$ and compute
\begin{align*}
    (n_{i1},\dots, n_{in})^{tr} &= E_i \cdot (a_1, \dots, a_n)^{tr} \quad \text{and} \\
    (n'_{i1},\dots, n'_{i,n'-n})^{tr} &= E'_i \cdot (a_{n+1}, \dots, a_{n'})^{tr}.
\end{align*}

\item[(5)] Return the bilinear maps $f_i \colon N_{1i} \times N_{2i} \longrightarrow M_i$
    for $i=1,\dots,\ell$, where $N_{1i}$ is the subgroup of~$N_1$ generated by $\{n_{i1}, \dots, n_{in}\}$,
    $N_{2i}$ is the subgroup of~$N_2$ generated by $\{n'_{i1}, \dots, n'_{in'}\}$,
    and $M_i$ is the subgroup of~$M$ generated by $f(n_{ij}, n'_{ik})$ for $j= 1,\dots,n$ and $k=1, \dots, n'-n$.

\end{enumerate}
\end{corollary}

Notice that all steps in this computation can be carried out in probabilistic polynomial time,
except for one integer prime factorization inside the algorithm of Step~(3).

\bigbreak
%
%

\section{Computing Decompositions of Finite $\ZZ$-Algebras}
\label{decomp_general}

In this section we combine the previous results and use the idempotents of the maximal
ring of scalars of a $\ZZ$-algebra $R$ to obtain a decomposition of $R$.
Recall that we denote the underlying $\ZZ$-module of the $\ZZ$-algebra~$R$ by~$R^+$.

\begin{definition}
Let $R_1,\dots,R_n$ be $\ZZ$-algebras.
\begin{enumerate}
\item[(a)] The {\bf direct product} of the $\ZZ$-algebras $R_1,\dots,R_n$ 
is the direct product $R_1^+ \times \cdots \times R_n^+$ of the underlying $\ZZ$-modules 
equipped with componentwise multiplication. It is denoted by $R_1 \times \cdots \times R_n$.

\item[(b)] A $\ZZ$-algebra is called \textbf{directly indecomposable} if it cannot 
be written as the direct product of two non-trivial $\ZZ$-subalgebras.

\end{enumerate}
\end{definition}

Let $R$ be a finite $\ZZ$-algebra. Our goal is to obtain a decomposition 
of~$R$ into a direct product of directly indecomposable $\ZZ$-subalgebras of~$R$. 
Let $R^2$ be the subgroup of $R^+$ generated by all products $a\cdot b$ with $a,b \in R^+$. 
Recall from Section~\ref{sec:max_ring_Zalg} that the maximal
ring of scalars $\mathfrak{S}(R)$ of~$R$ is a subring of the maximal ring of scalars 
of the bilinear map
$$
f_R^* \colon R^+/\Ann(R) \times R^+/\Ann(R) \longrightarrow R^2
$$
which satisfies the additional condition that the canonical map $R^2 \rightarrow R^+/\Ann(R)$ is
$\mathfrak{S}(R)$-linear. Using the primitive idempotents of~$\mathfrak{S}(R)$,
computed as in Section~\ref{idem},
we can now decompose the $\ZZ$-algebra $R/\Ann(R)$ as follows.

\begin{theorem}\label{decompose_R}{\textbf{(Direct Decomposition of Finite
$\ZZ$-Algebras Modulo Annihilator)}}
Let~$R$ be a finite $\ZZ$-algebra, and assume that we are given a complete set of primitive
idempotents $\{e_1, \dots, e_k\}$ of the maximal ring of scalars~$\mathfrak{S}(R)$.
For each $i \in \{1,\dots,k\}$, let $R_i = e_i\, (R/\Ann(R))$. 
Then we obtain a decomposition
$$
R/\Ann(R) \;=\; R_1 \times \cdots \times R_k
$$
of $R/\Ann(R)$ into a direct product of $\ZZ$-subalgebras.
\end{theorem}

\begin{proof}
Since $R^+/\Ann(R)$ is a faithful $\mathfrak{S}(R)$-module the primitive idempotents yield
a decomposition of the $\ZZ$-module $R^+/\Ann(R)$ into a direct sum of the 
$\ZZ$-submodules $R_i$. Therefore the map
$$
\phi \colon R_1 \times \cdots \times R_k \rightarrow R/\Ann(R), \quad 
(r_1, \dots, r_k) \mapsto r_1 + \cdots + r_k
$$
is a $\ZZ$-module isomorphism. Let us now show that the $\ZZ$-submodules $R_i$ 
are closed under multiplication.
We denote the canonical image of an element $x \in R$ in $R/\Ann(R)$ by $\bar{x}$.
Let $x,y \in R$ such that $\bar{x}, \bar{y} \in R_i$. 
Then these elements satisfy $e_i \bar{x} = \bar{x}$
and $e_i \bar{y} = \bar{y}$. Since the canonical homomorphism $\pi:R^2 \rightarrow R^+/\Ann(R)$
is $\mathfrak{S}(R)$-linear by the definition of $\mathfrak{S}(R)$, we get
$$
\bar{x}\, \bar{y} = \pi(xy) = \pi(f_R^*(\bar{x}, \bar{y})) = 
\pi(f_R^*(e_i \bar{x}, e_i \bar{y})) = \pi(e_i x y) = e_i \bar{x}\bar{y}.
$$
This shows that each $R_i$ is a subalgebra of $R/\Ann(R)$. We proceed to show that $\phi$ is a
$\ZZ$-algebra isomorphism. For $x,y \in R/\Ann(R)$ and $i\ne j$, we have
$$
e_ix \cdot e_j y = \pi(f_R^*(e_ix, e_jy))= e_i e_j \pi(f_R^*(x,y))=0.
$$
For $i=1, \dots, k$, let $r_i, s_i \in R_i$. Then we have $r_i = e_i r_i$ and $s_i = e_i s_i$, 
and the above equation implies
\begin{align*}
\phi((r_1,\dots,r_k))\, \phi((s_1, \dots, s_k)) 
&= (\tsum_{i=1}^k e_i r_i)(\tsum_{i=1}^k e_i s_i) = \\
&= \tsum_{i=1}^k e_i s_i r_i = \phi((r_1, \dots, r_k)(s_1, \dots, s_k)). \qedhere
\end{align*}
\end{proof}

Combining this theorem with the results of the proceeding sections, we get the following algorithm.
We assume that a $\ZZ$-algebra $R$ is given as in Remark~\ref{remark:algebra_input}. 
This means that~$R$ is generated by elements $a_1, \dots, a_n$, and we are given a presentation 
of the $\ZZ$-module~$R^+$. Furthermore, the multiplication
is represented by structure constants $s_{ijk}\in \ZZ$ for $i,j,k=1,\dots,n$ such that
$a_i a_j = \sum_{k=1}^m s_{ijk} a_k$.

\begin{corollary}\label{alg_decomposition}
Let $R$ be a finite $\ZZ$-algebra.
Assume that~$R$ is given as in Remark~\ref{remark:algebra_input}.
Then the following steps define an algorithm which computes a direct decomposition
$R/\Ann(R) = R_1 \times \cdots \times R_\ell$ into $\ZZ$-subalgebras of $R/\Ann(R)$.
\begin{enumerate}
\item[(1)] Use Corollary~\ref{Cor:maxring_algebra} to compute tuples
$(d_{ij}^{(1)}),\dots,(d_{ij}^{(t)}) \in (\ZZ^{2n})^{2n}$
corresponding to elements of~$\End(N_1) \times \End(N_2)$ which generate the
maximal ring of scalars~$\mathfrak{S}(f)$, together with an ideal
$I \subseteq \ZZ[x_1,\dots, x_t]$ such that $\mathfrak{S}(R) = \ZZ[x_1,\dots, x_t]/I$.
For $k=1, \dots, t$ let $A_k \in \Mat_{n,n}(\ZZ)$
be the matrix given by $(d_{i,j}^{(k)})_{i,j=1, \dots, n}$ and let
${A'}_k \in \Mat_{n,n}(\ZZ)$ be the matrix given by
$(d_{i,j}^{(k)})_{i,j=n+1, \dots, 2n}$.

\item[(2)] Use Proposition~\ref{alg_idempotents} to compute a complete set
$\{e_1, \dots, e_\ell\}$ of primitive idempotents of $\ZZ[x_1,\dots, x_t]/I$.

\item[(3)] For $i=1, \dots, \ell$, let $E_i = e_i(A_1, \dots, A_t)$ and
$E'_i = e_i(A'_1, \dots, A'_t)$ and compute
$$
r_{ij} \;=\; \eta^{-1}(E_i \cdot (0, \dots,b_j, \dots, 0)^{tr}, E'_i \cdot (0, \dots,c_j, \dots, 0)^{tr})
$$
for $j=1, \dots, n$, where $\eta : R \rightarrow R/\Ann_\lambda(R) \times R/\Ann_\rho(R)$ is the canonical group
homomorphism.

\item[(4)] Return the $\ZZ$-algebras $R_1, \dots, R_\ell$, where $R_i$ is the
subalgebra of $R/\Ann(R)$ generated by the residue classes of $r_{i1}, \dots, r_{in}$.

\end{enumerate}
\end{corollary}

Again, except for one integer prime factorization inside the algorithm of Step~(2), 
all computations can be performed in probabilistic polynomial time. Note that, in general, 
we cannot lift a direct decomposition of~$R/\Ann(R)$ to a direct decomposition of~$R$. 
Also, the direct decomposition of $R/\Ann(R)$ that we obtain
need not be directly indecomposable. Examples which show this are given below.

A class of examples for finite $\ZZ$-algebras for which the left and right annihilators coincide,
but which do not contain an identity element, is provided by finite dimensional Lie rings.

\begin{definition}
A \textbf{Lie ring} $L$ is a $\ZZ$-module together with a $\ZZ$-bilinear map
$[\ ,\ ] \colon L \times L \longrightarrow L$ such that
\begin{enumerate}
\item[(a)] $[x,x] = 0$ and

\item[(b)] $[x,[y,z]]+[z,[x,y]]+[y,[z,x]] = 0$
\end{enumerate}
for all $x,y,z \in L$. The bilinear map is called \textbf{Lie bracket}. A Lie ring is
called \textbf{finite dimensional} if it is finitely generated as an abelian group.
\end{definition}

Since we have $[x,y] = -[y,x]$ for all $x,y \in R$, a Lie ring~$R$, or more generally any
(anti-)commutative $\ZZ$-algebra, satisfies $\Ann_{\lambda}(R) = \Ann_\rho(R)$.
In the following example we apply Corollary~\ref{alg_decomposition}
to a concrete case.

\begin{example}
Let us consider the finite dimensional Lie ring~$L$ generated by $x_1,\dots,x_5$
with relations $3x_1 = 6x_2 = 3x_3 = 7x_4 = 7x_5 = 0$ whose multiplication is given by
\begin{alignat*}{3}
       & [x_1, x_2] = x_2,      \quad && [x_1, x_3] = 2x_5,\quad     && \\
       & [x_2, x_1] = 5x_2,     \quad && [x_2, x_3] = 3x_5,          && \\
       & [x_3, x_1] = 5x_5,     \quad && [x_3, x_2] = 4x_5,\quad     && [x_3, x_4] = x_3,\\
       & [x_4, x_3] = 2x_3
\end{alignat*}
and $[x_i,x_j] = 0$ in all remaining cases. When we compute the annihilators of~$L$, we obtain
$\Ann_\lambda(L) = \Ann_\rho(L) = \langle x_5, 3x_4 \rangle$. Following the steps of the algorithm in
Corollary~\ref{alg_decomposition} we get that the maximal ring of
scalars~$\mathfrak{S}(L)$ is generated by the endomorphisms $\phi_1, \phi_2 \in \End(L^+/\Ann(L))$
given by the matrices
\begin{equation*}
M_{\phi_1} = \begin{pmatrix}
           1 & 0 & 0 & 0 & 0 \\
           0 & 1 & 0 & 0 & 0 \\
           0 & 0 & 0 & 0 & 0 \\
           0 & 0 & 0 & 0 & 0 \\
           0 & 0 & 0 & 0 & 0
        \end{pmatrix}
\quad\hbox{and}\quad M_{\phi_2} =
\begin{pmatrix}
          0 & 0 & 0 & 0 & 0 \\
          0 & 0 & 0 & 0 & 0 \\
          0 & 0 & 1 & 0 & 0 \\
          0 & 0 & 0 & 0 & 0 \\
          0 & 0 & 0 & 0 & 0
\end{pmatrix}
\end{equation*}
These yield a presentation
\begin{equation*}
\mathfrak{S}(L) \cong \ZZ[y_1, y_2]/\langle 6, 3y_2, y_1 + y_2 - 1, y_2^2 - y_2 \rangle.
\end{equation*}
Using this presentation we compute the primitive idempotents $3$, $y_2$, and $2y_2+4$.
Finally, we get the direct decomposition $L/\Ann(L) = \langle 3\bar{x}_2 \rangle \times
\langle \bar{x}_3 \rangle \times \langle \bar{x}_1, 4\bar{x}_2 \rangle$.
\end{example}

The next two lemmas will come in handy when we analyze under the which conditions
the factors obtained in Corollary~\ref{alg_decomposition} are directly indecomposable.
In the following we consider $\ZZ$-algebras $R$ which satisfy $\Ann_\lambda(R) = \Ann_\rho(R)$.
Note that in this case we have $f^*_R = f_R$.

\begin{lemma}\label{lemma_in_ann}
Let $R$ be a finite $\ZZ$-algebra and let $S_1, S_2$ be subalgebras of $R/\Ann(R)$ 
such that~$R/\Ann(R)$ is the direct product of $S_1$ and $S_2$. Let $s_1, s_2 \in R$ 
be chosen such that $\bar{s}_1 \in S_1$ and $\bar{s}_2 \in S_2$.
Then we have $s_1 s_2, s_2 s_1\in \Ann(R) \cap R^2$.
\end{lemma}

\begin{proof}
Clearly, we have $s_1 s_2 , s_2 s_1 \in R^2$, and it remains to show $s_1 s_2, s_2 s_1 \in \Ann(R)$. 
By assumption, the $\ZZ$-linear map $\phi \colon S_1 \times S_2 \longrightarrow R/\Ann(R)$ 
given by $(a,b) \mapsto a+b$ is a $\ZZ$-algebra isomorphism. Therefore we get
$$
\phi((\bar{s}_1,0)(\bar{s}_1,\bar{s}_2)) = \phi((\bar{s}_1^2, 0)) =
    \bar{s}_1^2 = \bar{s}_1^2 + \bar{s}_1 \bar{s}_2 = \bar{s}_1 (\bar{s}_1+\bar{s}_2) =
    \phi((\bar{s}_1,0))\phi((\bar{s}_1,\bar{s}_2))
$$
This shows $\bar{s}_1\bar{s}_2 = 0$, and hence $s_1 s_2 \in \Ann(R)$. Analogously, we
obtain $s_2 s_1 \in \Ann(R)$.
\end{proof}

Under the following assumption a direct decomposition of the bilinear map $f_R$ 
yields a decomposition of~$\mathfrak{S}(R)$.

\begin{lemma}\label{lemma:max_ring_correspondence}
Let $R$ be a finite $\ZZ$-algebra that satisfies $\Ann_\lambda(R) = \Ann_\rho(R)$
and $\Ann(R) \cap R^2 = \{0\}$. Assume that $S_1, S_2$ are subalgebras of $R/\Ann(R)$
such that the bilinear map
$$
f_R :\; R^+/\Ann(R) \times R^+/\Ann(R) \longrightarrow R^2
$$
can be decomposed into a direct product of the non-trivial bilinear maps 
$f_{S_1} \times f_{S_2}$, where $f_{S_i} :\; S_i^+ \times S_i^+ \longrightarrow 
\langle f_R(S_i^+, S_i^+) \rangle$ is given by the restriction of~$f_R$ to
$S_i^+ \times S_i^+$. Then $\mathfrak{S}(R)$ can be decomposed into a direct 
product of non-trivial subrings.
\end{lemma}

\begin{proof}
By assumption, we have $R^+/\Ann(R) = S_1^+ \times S_2^+$. 
We denote the projection from $R^+/\Ann(R)$ to $S_i^+$ by $\pi_i \in \End(R^+/\Ann(R))$. 
Then, in view of the proof of \cite{Mya1}, Proposition~3.1, 
we have $\pi_i \in \mathfrak{S}(f_R)$, and $\pi_1,\pi_2$ are orthogonal idempotents 
in $\mathfrak{S}(f_R)$ with $1 = \pi_1+\pi_2$. 

Now we proceed to show that $\pi_i \in \mathfrak{S}(R) \subseteq \mathfrak{S}(f_R)$. 
Let $\psi : R^2 \rightarrow R^+/\Ann(R)$ be the canonical homomorphism. 
We want to show that $\psi(\pi_i(x)) = \pi_i(\psi(x))$ for all $x \in R^2$. 
By assumption, the group $R^2$ is generated by elements of the form 
$f_R(\bar{s}_1 + \bar{s}_2, \bar{s}'_1 + \bar{s}'_2)$ with $s_1, s'_1, s_2, s'_2 \in R$ 
such that $\bar{s}_1, \bar{s}'_1 \in S_1$ and $\bar{s}_2, \bar{s}'_2 \in S_2$. Then we have
$$
\psi(\pi_1(f_R(\bar{s}_1 + \bar{s}_2, \bar{s}'_1 + \bar{s}'_2))) =
\psi(f_R(\pi_1(\bar{s}_1 + \bar{s}_2), \pi_1(\bar{s}'_1 + \bar{s}'_2))) =
\psi(f_R(\bar{s}_1, \bar{s}'_1)) = \bar{s}_1 \bar{s}'_1.
$$
From Lemma~\ref{lemma_in_ann} we get $s_1s'_2+s_2s'_1 \in \Ann(R) \cap R^2 = \{0\}$, 
and therefore
$$
\pi_1(\psi(f_R(\bar{s}_1 + \bar{s}_2, \bar{s}'_1 + \bar{s}'_2))) =
\pi_1(\psi(s_1s'_1 + s_2s'_2)) ) = \pi_1(\bar{s}_1 \bar{s}'_1 + \bar{s}_2 \bar{s}'_2) =
\bar{s}_1 \bar{s}'_1.
$$
This proves $\pi_1 \in \mathfrak{S}(R)$, and $\pi_2 \in \mathfrak{S}(R)$ follows analogously.
Altogether, we see that $\mathfrak{S}(R) = \pi_1 \mathfrak{S}(R) \times \pi_2 \mathfrak{S}(R)$.
\end{proof}

The direct decomposition of $R/\Ann(R)$ obtained in Proposition~\ref{decompose_R} is, in general, 
not directly indecomposable, as we shall see below. However, under the following additional
hypothesis, it is.

\begin{proposition}\label{prop:indecomp_factors}
Let $R$ be a finite $\ZZ$-algebra with $\Ann_\lambda(R) = \Ann_\rho(R)$, and let
$R/\Ann(R) = R_1 \times \cdots \times R_k$ be the decomposition into $\ZZ$-subalgebras
obtained from the primitive idempotents $\{e_1, \dots, e_k \}$ of $\mathfrak{S}(R)$ as in
Proposition~\ref{decompose_R}. If $R^2 \cap \Ann(R) = \{ 0 \}$, then the factors $R_1,\dots,R_k$
are directly indecomposable.
\end{proposition}

\begin{proof}
Recall that $R_i = e_i(R/\Ann(R))$, and that the idempotents yield a decomposition 
of the bilinear map $f_R : R^+/\Ann(R) \times R^+/\Ann(R) \longrightarrow R^2$ 
into the direct product of the bilinear maps
$$
f_{R_i}: R_i^+ \times R_i^+ \longrightarrow \langle f_R(R_i^+, R_i^+) \rangle,
\quad (e_i \bar{x}, e_i \bar{y}) \mapsto e_i x y.
$$
Assume that $R_i$ is the direct product $R_i = S_1 \times S_2$ of two non-trivial subalgebras
$S_1,S_2 \subseteq R_i$. For $j=1,2$, let~$M_j$ be the $\ZZ$-submodule 
generated by $f_R(S_j^+,S_j^+)$. We show that the bilinear map $f_{R_i}$ can
be further decomposed into the direct product of the two non-trivial bilinear maps
$ f_{S_j}:\; S_j^+ \times S_j^+ \longrightarrow M_j $ given by the restriction 
of~$f_R$ to $S_j^+ \times S_j^+$.

By assumption, we have $R_i^+ = S_1^+ \oplus S_2^+$.
Next we show that $\langle f_R(R_i^+, R_i^+) \rangle = M_1 \oplus M_2$. 
Let $\pi : R^2 \rightarrow R^+/\Ann(R)$ be the canonical homomorphism. 
Since $S_1 \cap S_2 = 0$ and $\pi(M_j) \subseteq S_j$, we have
$\pi(M_1 \cap M_2) = 0$, and hence
$$
M_1 \cap M_2 \subset \Ann(R) \cap R^2 = \{ 0 \}.
$$
To show that $\langle f_R(R_i^+, R_i^+) \rangle = M_1 + M_2$, let $\bar{s}_1 + \bar{s}_2$ and
$\bar{s}'_1 + \bar{s}_2'$ be two elements of~$R_i$ with $\bar{s}_1, \bar{s}'_1 \in S_1$ and
$\bar{s}_2, \bar{s}'_2 \in S_2$. Then we have
\begin{align*}
    f_R(\bar{s}_1 + \bar{s}_2, \bar{s}_1' + \bar{s}_2')
    &= f_R(\bar{s}_1, \bar{s}_1') + f_R(\bar{s}_2, \bar{s}_2') + s_1 s'_2 + s_2 s'_1 = \\
    &= f_{S_1}(\bar{s}_1, \bar{s}_1') + f_{S_2}(\bar{s}_2, \bar{s}_2')
\end{align*}
since $s_1 s'_2 + s_2 s'_1 \in \Ann(R) \cap R^2 = \{0\}$ by Lemma~\ref{lemma_in_ann}. 
This proves $\langle f_R(R_i^+, R_i^+) \rangle = M_1 \oplus M_2$, and it also shows that
$f_{R_i}$ is the direct product of the bilinear maps $f_{S_1}$ and $f_{S_2}$.
By Lemma~\ref{lemma:max_ring_correspondence} the ring $e_i \mathfrak{S}(R)$ can then be decomposed
into a direct product of subrings. Therefore the primitive idempotent~$e_i$ is the sum of two
non-trivial orthogonal idempotents, which is a contradiction.
\end{proof}

Our next example shows that the direct decomposition of $R/\Ann(R)$ given in
Proposition~\ref{decompose_R} is, in general, not directly indecomposable.

\begin{example}\label{ex:not_a_lifting}
Let us consider the finite dimensional Lie ring~$R$ generated by $x_1,\dots,x_5$
with multiplication given by
$[x_1, x_2] = x_5$, $[x_3, x_4] = x_5$, and $[x_i, x_j] = 0$ in all remaining cases.
This implies $\Ann(R) = \langle x_5 \rangle$ and $R^2 = \langle x_5 \rangle$.
Clearly, the product $R/\Ann(R) = \langle \bar{x}_1, \bar{x}_2 \rangle
\times \langle \bar{x}_3, \bar{x}_4 \rangle$ is a direct decomposition of $R/\Ann(R)$.
However, the maximal ring of scalars~$\mathfrak{S}(R)$ is indecomposable, since the
$\ZZ$-module $R^2$ is directly indecomposable. Therefore we obtain no decomposition
of $R/\Ann(R)$ using Proposition~\ref{decompose_R}.
\end{example}

Moreover, in general, we cannot lift a direct decomposition of~$R/\Ann(R)$ to a direct
decomposition of~$R$, as the following two examples show.

\begin{example}\label{ex:not_liftable1}
Let $R$ be the free $\ZZ$-module generated by $x_1,\dots,x_6$. We define a commutative
multiplication on~$R$ by
$$
x_1 x_3 = x_4,\quad x_2 x_3 = x_4,\quad x_3 x_4 = x_3, \quad x_5 x_6 = x_1 - x_2
$$
and $x_i x_j = 0$ in all other cases. Then we obtain $\Ann(R) = \langle x_1 - x_2 \rangle$,
since we have $(x_1 - x_2) x_3 = x_4 - x_4 = 0$.
The maximal ring of scalars~$\mathfrak{S}(R)$ is generated by the endomorphisms
$\phi_1, \phi_2 \in \End(R)$ given by the matrices
\begin{equation*}
M_{\phi_1} =  \begin{pmatrix}
        1 & 0 & 0 & 0 & 0 & 0 \\
        1 & 0 & 0 & 0 & 0 & 0 \\
        0 & 0 & 1 & 0 & 0 & 0 \\
        0 & 0 & 0 & 1 & 0 & 0 \\
        0 & 0 & 0 & 0 & 0 & 0 \\
        0 & 0 & 0 & 0 & 0 & 0
\end{pmatrix}
\quad\hbox{\rm and}\quad
M_{\phi_2} =  \begin{pmatrix}
        0 & 0 & 0 & 0 & 0 & 0 \\
        0 & 0 & 0 & 0 & 0 & 0 \\
        0 & 0 & 0 & 0 & 0 & 0 \\
        0 & 0 & 0 & 0 & 0 & 0 \\
        0 & 0 & 0 & 0 & 1 & 0 \\
        0 & 0 & 0 & 0 & 0 & 1
\end{pmatrix}
\end{equation*}
Then the maximal ring of scalars has the presentation
\begin{equation*}
\mathfrak{S}(R) \cong \ZZ[y_1, y_2] / \langle y_1^2-y_1, y_1y_2, y_2^2-y_2, y_1+y_2-1 \rangle,
\end{equation*}
and from this we compute the primitive idempotents $-y_2+1$ and $y_2$.

Finally, we get the direct decomposition
$R/\Ann(R) = \langle \bar{x}_1, \bar{x}_2, \bar{x}_3, \bar{x}_4 \rangle \times 
\langle \bar{x}_5,  \bar{x}_6 \rangle$.
This decomposition cannot be lifted to a direct decomposition of~$R$ because of
$x_1 - x_2 \in \langle x_5, x_6 \rangle \cap \langle x_1, x_2, x_3, x_4 \rangle$.
\end{example}

In the final example of this section we see a case in which a decomposition of $R/\Ann(R)$
cannot be lifted to~$R$ for a different reason.

\begin{example}\label{ex:not_liftable2}
Let $R$ be the free $\ZZ$-module generated by $x_1,\dots,x_5$. We define a commutative
multiplication on~$R$ by letting
$$
x_1^2 = 2x_4, \quad x_1 x_2 = x_2,\quad x_2^2 = x_2, \quad x_3^2 = x_3
$$
and $x_i x_j = 0$ in all other cases. Then we have $\Ann(R) = \langle x_4, x_5 \rangle$
and $R^2 = \langle x_2, x_3, 2x_4 \rangle$. The maximal ring of scalars~$\mathfrak{S}(R)$
is generated by the endomorphisms $\phi_1, \phi_2 \in \End(R)$ given by the matrices
\begin{equation*}
M_{\phi_1} =   \begin{pmatrix}
       1 & 0 & 0 & 0 & 0 \\
       0 & 1 & 0 & 0 & 0 \\
       0 & 0 & 0 & 0 & 0 \\
       0 & 0 & 0 & 0 & 0 \\
       0 & 0 & 0 & 0 & 0
\end{pmatrix}
\quad\hbox{\rm and}\quad
M_{\phi_2} =  \begin{pmatrix}
       0 & 0 & 0 & 0 & 0 \\
       0 & 0 & 0 & 0 & 0 \\
       0 & 0 & 1 & 0 & 0 \\
       0 & 0 & 0 & 0 & 0 \\
       0 & 0 & 0 & 0 & 0
\end{pmatrix}
\end{equation*}
Here we calculate the presentation
\begin{equation*}
\mathfrak{S}(R) \cong \ZZ[y_1, y_2]/\langle y_1^2 - y_1, y_1y_2, y_2^2 - y_2, y_1 + y_2 - 1 \rangle,
\end{equation*}
and from this we determine the primitive idempotents $-y_2+1$ and $y_2$ of~$\mathfrak{S}(R)$.

Finally, we obtain the direct decomposition 
$R/\Ann(R) = \langle \bar{x}_1, \bar{x}_2 \rangle \times
\langle \bar{x}_3 \rangle$. To lift this decomposition to a decomposition of~$R$, 
we would need to find a direct complement of $\Ann(R) \cap R^2 = \langle 2x_4 \rangle$ in 
$\Ann(R) = \langle x_4, x_5 \rangle$. However, it is clear that $\langle 2x_4 \rangle$ 
is not a direct summand in $\langle x_4, x_5 \rangle$.
\end{example}

In view of these examples we can lift the decompositions of~$R/\Ann(R)$ to
decompositions of~$R$ only under suitable additional assumptions. One setting where
one could look for such hypotheses are Lie rings associated to finitely generated 
nilpotent groups.

Let us denote the members of the lower central series
of a group~$G$ by $\gamma_i(G)$, i.e., we let $\gamma_1(G) = G$
and $\gamma_{k+1}(G) = [\gamma_k(G), G]$ for $k\ge 1$. The group~$G$ is called {\bf nilpotent}
if $\gamma_n(G) = \{e\}$ for some $n\ge 0$. The smallest such~$n$ is then called
the {\bf nilpotency class} of~$G$. Let us also recall the following construction.

\begin{definition}
Let $G$ be a nilpotent group. The additive group of the \textbf{associated Lie ring} $L(G)$
is given by the direct sum $L(G) = \bigoplus_{k=1}^{\infty} \gamma_k/\gamma_{k+1}$.
The direct summand $\gamma_k/\gamma_{k+1}$
is called the \textbf{homogeneous component of weight}~$k$ of~$L(G)$.

Multiplication of elements of homogeneous components of~$L(G)$ is defined by
$$
[a+\gamma_{i+1}, b+\gamma_{j+1}] = [a,b] + \gamma_{i+j+1}
$$
for $a\in \gamma_i(G)$ and $b\in \gamma_j(G)$.
Then this multiplication is extended to~$L(G)$ by linearity.
\end{definition}

Let us consider the case of nilpotent groups of class~2, i.e., groups~$G$ 
such that $[G, [G,G]] = \{e\}$.

\begin{remark}
If $G$ is nilpotent of class~2, the definition of $L(G)$ simplifies 
to $L(G) = G/[G,G] \times [G,G]$. An element $a[G,G]+b \in L(G)$ is in 
$\Ann(L(G))$ if and only if
$[a\, [G,G]+b,\, c\, [G,G]+d] = [a\,[G,G],\, c\,[G,G]] = [a,c] = e$ for all $c\in G$.
This is equivalent to $a \in \Cen(G)$, where $\Cen(G)$ denotes the center of~$G$.
Hence we have
$$
\Ann(L(G)) = \Cen(G)/[G,G] \times [G,G]
$$
and the bilinear map associated to $L(G)$ is given by
$$
f_G : G/\Cen(G) \times G/\Cen(G) \longrightarrow [G,G], \quad
(g\Cen(G),h\Cen(G)) \mapsto [g,h].
$$
Note that the canonical group homomorphism $[G,G] \rightarrow G/\Cen(G)$ is trivial,
and therefore $\mathfrak{S}(f_G)$ linear. This implies $\mathfrak{S}(L(G)) = \mathfrak{S}(f_G)$ 
by Theorem~\ref{max_ring_algebra_existence}. Now the primitive idempotents 
of $\mathfrak{S}(L(G))$ yield a decomposition of~$f_G$ into a direct product of 
indecomposable bilinear maps. In particular, we obtain a decomposition of the 
abelian group $G/\Cen(G)$ into cyclic subgroups. 

At this point a number of further questions arise.
Can we deduce a decomposition of~$G$ from the decomposition of the bilinear map $f_G$? 
Can one find decompositions for nilpotent groups of higher nilpotency classes using this method?
We leave this study to future research.
\end{remark}

\medskip\noindent
{\bf Acknowledgements.} The first and third author thank the Stevens Institute,
Hoboken, USA, for its hospitality and support during part of the preparation 
of this paper.

\bigbreak

\end{document}